\documentclass[leqno,draft,11pt]{article}
\usepackage{amssymb,amsmath,latexsym,theorem,a4wide}

\allowdisplaybreaks[2]

\newcommand\M{\mathit}

\newcommand\eq{\leftrightarrow}
\newcommand\LOR{\bigvee}
\newcommand\ET{\bigwedge}
\newcommand\model{\vDash}
\newcommand\bool{\mathcal B}
\newcommand\mbool{\mathcal M}

\newcommand\fii{\varphi}
\newcommand\tet{\theta}

\newcommand\p[1]{\langle#1\rangle}
\newcommand\lh[1]{\lvert#1\rvert}
\newcommand\sset{\subseteq}
\newcommand\Sset{\supseteq}
\newcommand\res{\mathbin\restriction}

\newcommand\fl[1]{\lfloor#1\rfloor}
\newcommand\cl[1]{\lceil#1\rceil}
\newcommand\fdiv{\genfrac\lfloor\rfloor{}{}}
\newcommand\dive[2]{\fdiv{#1}{2^{#2}}}
\newcommand\tdive[2]{\fl{#1/2^{#2}}}
\newcommand\half[1]{\fl{#1/2}}
\newcommand\str[1]{[\mkern-2.7mu[#1]\mkern-2.7mu]}
\newcommand\rsuv{\M{RSUV}}
\providecommand\dotminus{\mathbin{\mathchoice
       {\scriptstyle\dodotminus\displaystyle}%
       {\scriptstyle\dodotminus\textstyle}%
       {\scriptscriptstyle\dodotminus\scriptstyle}%
       {\scriptscriptstyle\dodotminus\scriptscriptstyle}}}
\def\dodotminus#1{\dot{\smash{#1-}}}

\newcommand\seq{\Longrightarrow}

\newcommand\ru{\mathrel/}

\DeclareMathOperator\Th{Th}

\newcommand\cxt{\mathrm}
\newcommand\ptime{\cxt P}
\newcommand\ph{\cxt{PH}}
\newcommand\fp{\cxt{FP}}
\newcommand\nc{\cxt{NC}}
\newcommand\nci{\nc^1}
\newcommand\tc{\cxt{TC}^0}

\newcommand\sig{\Sigma^b_}
\newcommand\pii{\Pi^b_}
\newcommand\delt{\Delta^b_}
\newcommand\Sig{\Sigma^B_}
\newcommand\st{\expandafter\hat}
\newcommand\sih{\st\sig}
\newcommand\pih{\st\pii}
\newcommand\sigp{\Sigma_1^*}
\newcommand\piip{\Pi_1^*}
\newcommand\idz{I\Delta_0}
\newcommand\bas{\M{BASIC}}
\newcommand\btc{\M{BTC}^0}
\newcommand\pv{PV}
\newcommand\tpv{PV_1}
\newcommand\vtc{\M{VTC}^0}
\newcommand\BB{\M{BB}}
\newcommand\sch[1]{\ensuremath{\M{#1}}}
\newcommand\lhmin{\sch{LMIN}}
\newcommand\Min{\sch{MIN}}
\newcommand\lind{\sch{LIND}}
\newcommand\ind{\sch{IND}}
\newcommand\pind{\sch{PIND}}
\newcommand\indf{\sch{IND^-}}
\newcommand\pindf{\sch{PIND^-}}
\newcommand\indr{\sch{IND^R}}
\newcommand\pindr{\sch{PIND^R}}
\newcommand\ppind{\sch{(P)IND}}
\newcommand\ppindf{\sch{(P)IND^-}}
\newcommand\ppindr{\sch{(P)IND^R}}
\newcommand\llmin{\sch{(L)MIN}}
\newcommand\rfn{\mathrm{RFN}}

\newcommand\G{G}
\newcommand\siq{\Sigma^q_}
\newcommand\piq{\Pi^q_}

\newcommand\stm{\mathbb N}

\newcommand\ob{\overline}
\def\cput(#1)#2{\put(#1){\hbox to0pt{\hss#2\hss}}}

\newcommand\bme{\hskip.75em\relax}
\newcommand\noproof{\leavevmode\unskip\bme\vadjust{}\nobreak\hfill$\qed$\par}
\newcommand\qed{\Box}
\newenvironment{Pf}[1][]
  {\par\noindent\textit{Proof\optpar{#1}:}\bme\ignorespaces}
  {\noproof\pagebreak[2]\vskip\medskipamount\ignorespacesafterend}
\def\optpar#1{\ifx\relax#1\relax\else\ #1\fi}
\newcommand\qedhere{\relax\ifmmode\eqno\qed\expandafter\aftergroup
                   \else\noproof\fi\noqed}
\newcommand\noqed{\let\noproof\relax}

\theoremstyle{plain}
\newtheorem{Thm}{Theorem}[section]
\newtheorem{Prop}[Thm]{Proposition}
\newtheorem{Cor}[Thm]{Corollary}
\newtheorem{Lem}[Thm]{Lemma}
\newtheorem{Obs}[Thm]{Observation}
\newtheorem{Que}[Thm]{Question}
\newtheorem{Cl}{Claim}[Thm]

\theorembodyfont\upshape
\newtheorem{Def}[Thm]{Definition}
\newtheorem{Rem}[Thm]{Remark}
\newenvironment{Pf*}{\let\qed\qedCl\Pf}\endPf

\numberwithin{figure}{section}

\newcommand\allowhyphens{\nobreak\hskip0pt\relax}
\DeclareRobustCommand*\magiclparen{\ifmmode(\else\textup(\allowhyphens\fi}
\DeclareRobustCommand*\magicrparen{\ifmmode)\else\textup)\fi}
\newcommand\magicparon{\catcode`\(\active\catcode`\)\active}
\newcommand\magicparoff{\catcode`\(12 \catcode`\)12 }
\magicparon
\let (=\magiclparen  \let )=\magicrparen

\author{Emil Je\v r\'abek\dedic\\[\medskipamount]
The Czech Academy of Sciences, Institute of Mathematics\\
\small \v Zitn\'a 25,
115\:67 Praha 1,
Czech Republic,
email: \texttt{jerabek@math.cas.cz}
}

\newcommand\dedic{\thanks{Supported by Center of Excellence CE-ITI under the grant
P202/12/G061 of GA \v CR. The Institute of Mathematics of the Czech Academy of
    Sciences is supported by RVO: 67985840.}}

\title{Induction rules in bounded arithmetic}

\begin{document}
\maketitle

\begin{abstract}
We study variants of Buss's theories of bounded arithmetic axiomatized by induction schemes disallowing the use of
parameters, and closely related induction inference rules. We put particular emphasis on $\pih i$~induction schemes,
which were so far neglected in the literature. We present inclusions and conservation results between the systems
(including a witnessing theorem for $T^i_2$ and~$S^i_2$ of a new form),
results on numbers of instances of the axioms or rules, connections to reflection principles for quantified
propositional calculi, and separations between the systems.


\end{abstract}

\section{Introduction}\label{sec:introduction}

Commonly studied theories of arithmetic, weak and strong alike, are typically axiomatized by variants of induction or
other axiom schemes (comprehension, collection, \dots) restricted to suitable classes of formulas, where these
formulas may freely use \emph{parameters}: arbitrary numbers or other objects manipulated by the theory that enter the
induction formula by means of free variables, unrelated to the induction variable. This generally makes the theories
robust in their formal properties, and intuitive to work with. Nevertheless,
induction schemes without parameters proved fruitful to study in the context of strong subtheories of Peano arithmetic
($\Sigma_n$-induction), revealing a landscape of strange, and yet familiar systems: see e.g.\ Kaye, Paris, and
Dimitracopoulos~\cite{kpd}, Adamowicz and Bigorajska~\cite{adam-big:isig1-},
Bigorajska~\cite{bigo:ipi1-}, Beklemishev~\cite{bekl:indru,bekl:parfree}, and
Cord\'on-Franco and Lara-Mart\'\i n~\cite{cf-lm:loc-ind}.

On the one hand, the parameter-free
induction schemes $I\Sigma_n^-$ and $I\Pi_n^-$ (for $n\ge1$) are close to the original schemes with parameters $I\Sigma_n$, as
the theories are conservative over each other with respect to large classes of sentences (though the correspondence is
a bit off, as $I\Pi_{n+1}^-$ is on the same level as $I\Sigma_n$ and~$I\Sigma_n^-$). On the other hand, there are
substantial differences: as already alluded to, the $\Pi_n$~schemes without parameters become genuinely distinct from
(and weaker than) the matching $\Sigma_n$~schemes, whereas $I\Sigma_n=I\Pi_n$; neither
$I\Sigma_n^-$ nor~$I\Pi_n^-$ are finitely axiomatizable, in contrast to $I\Sigma_n$.

The parameter-free schemes $I\Sigma_n^-$ and~$I\Pi_n^-$ are intimately connected to induction \emph{rules}
$I\Sigma_n^R$ and~$I\Pi_n^R$: here, instead of theories generated just by axioms on top of the usual
rules of first-order logic, we consider a form of induction as an additional (Hilbert-style) rule of inference. It turns out
$I\Sigma_n^-$ is the weakest theory all of whose extensions are closed under~$I\Sigma_n^R$, and likewise for~$\Pi_n$. 
An important role in the analysis of $I\Sigma_n^-$ and~$I\Pi_n^-$ is played by \emph{reflection principles} for
fragments of arithmetic~\cite{bekl:indru,bekl:parfree}: while $I\Sigma_n$ is equivalent to a certain uniform (global)
reflection principle, the theories $I\Sigma_n^-$ and~$I\Pi_n^-$ can be characterized using relativized \emph{local}
reflection principles. There are also intricate connections relating the nesting of applications of rules and the
number of instances of axioms. As an alternative to reflection principles, parameter-free induction schemes can 
be analysed using \emph{local induction}~\cite{cf-lm:loc-ind}.

In contrast to all these results, much less is known about parameter-free induction axioms and induction rules in the
context of bounded arithmetic: the early work of Kaye~\cite{kaye:dio} introduced the parameter-free subtheories
$IE_i^-$ of~$\idz$, while the only investigation of parameter-free Buss's theories was done by Bloch~\cite{bloch}, who studied
proof-theoretically $\sig i$~parameter-free induction rules\footnote{Warning: the proof of Theorem~27, which
effectively claims that $\sih i\text-\ppindf\equiv\sih i\text-\ppindr$, is incorrect.} in a sequent formalism, and
Cord\'on-Franco, Fernand\'ez-Margarit, and Lara-Mart\'\i n~\cite{cffmlm}, whose main results concern conservativity of
the theories $S^i_2$ and~$T^i_2$ over the parameter-free and induction-rule versions of $\sih i$-\pind\ and~$\sih
i$-\ind, and conservativity of $\BB\sig i$ over its rule version. They rely on model-theoretic methods exploiting
variants of existentially closed models.

The purpose of this paper is to study parameter-free versions of Buss's theories in a more systematic way, filling in
various gaps in our knowledge to obtain a more complete picture. Some highlights are as follows. We will
investigate $\pih i$~schemes and rules, which were so far entirely ignored in the literature, alongside their $\sih
i$~counterparts; in particular, we will prove conservation results of $T^i_2$ and~$S^i_2$ over $\pih i$-\ppindf.
We try to get as complete a description of the relationships among the
systems in question as possible; to this end, we also include tentative separation results (conditional or
relativized). While bounded arithmetic is too weak to prove the consistency of interesting first-order
theories, it has a well-known connection to propositional proof systems; in accordance with this,
we will present characterizations of our systems in terms of variants of reflection principles for fragments of the
quantified propositional sequent calculus. We also include some results on the nesting of rules, namely conditions ensuring that closure under the
induction rules collapses to unnested closure, and conservation results of $n$~instances of parameter-free induction
axioms over $n$~applications of induction rules.

The paper is organized as follows. After some preliminary background in Section~\ref{sec:notation}, we introduce in
Section~\ref{sec:main-fragments} the main axioms and rules that we are interested in, and we prove some of their
elementary properties---primarily reductions between the rules (Theorem~\ref{prop:basicred}), but also a result on a
collapse of $\pih i$-\ppindr\ to unnested applications (Theorem~\ref{prop:singlepihindr}). We discuss various variants of
the axioms and rules in Section~\ref{sec:variants}, and we show them mostly equivalent to our main systems
(Proposition~\ref{prop:variants}).

The most substantial technical part of the paper comes in Section~\ref{sec:conserv}, which is devoted to conservation
results. We recall the conservation of $T^i_2$ and~$S^i_2$ over $\sih i$-\ppindr\ (Theorem~\ref{thm:conssig})
from~\cite{bloch,cffmlm}, and we set out to prove an analogous conservation result over $\pih i$-\ppindr\
(Theorem~\ref{thm:picons}). A key part of the proof is a new witnessing theorem for $\forall\exists\forall\sih{i-1}$
consequences (and $\forall\exists\forall\sih i$ consequences) of $T^i_2$ and~$S^i_2$, which may be of independent
interest (Theorem~\ref{thm:conservp} and Proposition~\ref{prop:conservs}). We obtain conservation results over $\Gamma$-\ppindf, summarized in
Corollary~\ref{cor:consindf}, and a result on collapse of nesting of $\sih i$-\ppindr\ (Theorem~\ref{thm:sihinst}). We also prove
more direct conservation results of $T+\Gamma\text-\ppindf$ over $T+\Gamma\text-\ppindr$ for arbitrary theories~$T$
(Theorem~\ref{thm:consindfr}).

We discuss connections to propositional proof systems in Section~\ref{sec:pps}, the main result being a
characterization of $\Gamma$-\ppindr\ and $\Gamma$-\ppindf\ in terms of reflection principles for quantified
propositional calculi (Theorem~\ref{thm:ppseq}). Section~\ref{sec:separations} is devoted to separations between our
systems: we present some conditional separations in Section~\ref{sec:unrel-separ}, and unconditional relativized
separations in Section~\ref{sec:relat-separ}. We conclude the paper with a few remarks in Section~\ref{sec:conclusion}.

\section{Notation and preliminaries}\label{sec:notation}

We assume the reader is familiar with the basics of bounded arithmetic. We will work in the framework of Buss's
one-sorted theories $S^i_2$ and~$T^i_2$, as presented e.g.\ in Buss~\cite{buss}, H\'ajek and Pudl\'ak~\cite[Ch.~V]{hp},
or Kraj\'\i\v cek~\cite{book}. It would not be too difficult to adapt our results to the setting of two-sorted theories
$V^i$ as in Cook and Nguyen~\cite{cook-ngu}, but we find the one-sorted setting simpler to use for the present purpose.

In order not to get bogged down in trivial technicalities, we will employ a robust base theory in a rich language in
place of Buss's $\bas$: let $\btc$ denote the basic first-order theory for~$\tc$, in a language $L_{\tc}$ with function
symbols for all $\tc$~functions so that $\btc$ is a universal theory. We are not very particular about its exact
definition; for example, we may axiomatize it as the theory $\delt1\text-\M{CR}$ of Johannsen and
Pollett~\cite{joh-pol:d1cr} expanded with function symbols for all $\sig1$-definable functions of the theory,
or as the equivalent theory $\M{TTC}^0$ of Clote and Takeuti~\cite{cl-tak:tc0}.
Note that $\btc$ is $\rsuv$-isomorphic to the theory $\vtc$ (or rather, $\overline{\M{VTC}}^0$) of Cook and
Nguyen~\cite{cook-ngu}. Unless stated otherwise, we will assume all first-order theories to be formulated in~$L_{\tc}$
and to extend~$\btc$.

If $\Gamma$ is a (possibly empty) set of sentences, and $\fii$ a sentence, we write $\Gamma\vdash\fii$ if $\fii$ is
provable in the theory $\btc+\Gamma$. We may omit outermost universal quantifiers when writing down $\Gamma$ or~$\fii$,
as is the customary fashion. We may also write $\Gamma\vdash\Delta$ for a set of sentences~$\Delta$, meaning
$\Gamma\vdash\fii$ for all $\fii\in\Delta$. We stress that $\btc+\Gamma$ is only closed under the
standard deduction rules of first-order logic (i.e., it includes logically valid sentences, and it is closed under
modus ponens); it is not supposed to be closed under the $\delt1\text-\M{CR}$ rule even if we define $\btc$ as
in~\cite{joh-pol:d1cr}.

Let $\sih i$ and~$\pih i$ denote the classes of \emph{strict} $\sig i$ and~$\pii i$ formulas in~$L_{\tc}$: that is,
$\sih0=\pih0=\sig0=\pii0$ is the class of sharply bounded formulas, and for $i>0$, a $\sih i$~formula ($\pih
i$~formula) consists of
$i$~alternating (possibly empty) blocks of bounded quantifiers followed by a $\sig0$~formula, where the first block is
existential (universal, resp.). Equivalently, we could further restrict the blocks to a single quantifier apiece. 
Note that every $\sig0$~formula is equivalent to an atomic formula in~$\btc$. The class of all bounded
formulas is denoted~$\sig\infty$.

We will combine notations such as $\sih i$ and~$\pih i$ with symbolic prefixes denoting \emph{unbounded} quantifiers: for
example, $\forall\exists\sih i$ denotes the class of formulas (in most contexts, sentences) consisting of a block of
universal quantifiers, followed by a block of existential quantifiers, followed by a $\sih i$~formula.

Let $\Gamma$ be a class of sentences, and $T$ a theory. The \emph{$\Gamma$-fragment of~$T$} is the theory axiomatized
by $\btc+\{\fii\in\Gamma:T\vdash\fii\}$. If $S$ is another theory, $T$ is \emph{$\Gamma$-conservative over~$S$} if the
$\Gamma$-fragment of~$T$ is included in~$S$.

Let $\sigp$ denote the least class of formulas that includes bounded formulas, and is closed under existential
and bounded universal quantifiers; $\piip$ denotes the dual class. A model-theoretic characterization of these classes
is that $\piip$~formulas are preserved downwards in cuts, and $\sigp$~formulas upwards.
\begin{Thm}[Parikh]\label{thm:parikh}
Let $T$ be a $\piip$-axiomatized extension of~$\btc$, and $\fii\in\sigp$. If $T\vdash\forall x\,\exists y\,\fii(x,y)$,
there exists a term~$t$ such that $T\vdash\forall x\,\exists y\le t(x)\,\fii(x,y)$.
\noproof\end{Thm}

We will occasionally use that $\sigp$-sentences true in the standard model of arithmetic~$\stm$ are provable in~$\btc$.

Another fundamental tool for studying systems of bounded arithmetic is Buss's witnessing theorem. We are actually not
interested in witnessing per se, but in the following consequence:
\begin{Thm}[Buss]\label{thm:buss}
For any $i\ge0$, $S^{i+1}_2$ is a $\forall\sih{i+1}$-conservative extension of~$T^i_2$.
\noproof\end{Thm}
We will in fact use it in an ostensibly stronger form:
\begin{Cor}\label{cor:buss}
For any $i\ge0$ and $T\sset\exists\forall\sih i$, $S^{i+1}_2+T$ is $\forall\exists\sih{i+1}$-conservative over $T^i_2+T$.
\end{Cor}
\begin{Pf}
Assume that $S^{i+1}_2+\exists u\,\forall z\,\psi(u,z)\vdash\forall x\,\exists y\,\fii(x,y)$, where $\psi\in\sih i$, and
$\fii\in\sih{i+1}$. Then $S^{i+1}_2$ proves $\forall x,u\,\exists y\,\bigl(\neg\psi(u,y)\lor\fii(x,y)\bigr)$. By
Parikh's theorem, we may bound the $y$~quantifier by a term in~$x$ and~u, which makes the statement (equivalent to) a
$\forall\sih{i+1}$~sentence. Thus, it is provable in~$T^i_2$ by Theorem~\ref{thm:buss}, and this implies $T^i_2+\exists
u\,\forall z\,\psi(u,z)\vdash\forall x\,\exists y\,\fii(x,y)$.
\end{Pf}

Our basic objects of study will be \emph{rules} rather than just axiom schemes. Here, a rule~$R$ is a set of
pairs $\p{\Gamma,\fii_0}$, where $\fii_0$ is a sentence, and $\Gamma=\{\fii_1,\dots,\fii_n\}$ is a finite set of
sentences; each $\p{\Gamma,\fii_0}\in R$ is called an \emph{instance} of~$R$, and will be written more
conspicuously as $\Gamma\ru\fii_0$, or
\begin{equation}\label{eq:22}
\frac{\fii_1\quad\fii_2\quad\dots\quad\fii_n}{\fii_0}.
\end{equation}
The instance above is \emph{$n$-ary}. We will identify axiom schemes with $0$-ary rules. Again, we will often omit
outermost universal quantifiers from the sentences~$\fii_i$ when writing down rules like~\eqref{eq:22}.

If $T$ is a theory, and $R$ a rule, then $T+R$ denotes the least theory~$T'$ (i.e., deductively closed set of
sentences) which includes $T$, and which is \emph{closed under~$R$}, meaning that for any instance $\Gamma\ru\fii$
of~$R$, if $\Gamma\sset T'$, then $\fii\in T'$. 
We may stratify this definition by counting the nesting depth of applications of the rules. Let $[T,R]$ denote the
closure of~$T$ under \emph{unnested} applications of $R$-instances, i.e., the theory axiomatized by
\[T\cup\bigl\{\fii:\Gamma\ru\fii\in R,T\vdash\Gamma\bigr\},\]
and we define $[T,R]_0=T$, $[T,R]_{n+1}=[[T,R]_n,R]$ by induction on~$n\in\omega$. Notice that $T+R=\bigcup_n[T,R]_n$.
See also Remark~\ref{rem:instances}.

A rule $R$ is \emph{weakly reducible} to a rule~$S$ if $T+R\sset T+S$ for all theories~$T$, and $R$ and~$S$ are
\emph{weakly equivalent} if they are weakly reducible to each other.
Likewise, $R$ is \emph{reducible} to~$S$, written $R\le S$, if $[T,R]\sset[T,S]$ for every theory~$T$, and $R$
and~$S$ are \emph{equivalent}, written $R\equiv S$, if $R\le S\le R$.
While this definition speaks about arbitrary theories~$T$, the only interesting cases are theories axiomatized by
premises of instances of the rules:
\begin{Obs}\label{obs:rule-red}
Let $R$ and~$S$ be rules.
\begin{enumerate}
\item\label{item:45}
$R$ is weakly reducible to~$S$ iff $\fii\in\btc+\Gamma+S$ for all instances $\Gamma\ru\fii$ of~$R$.
\item\label{item:46}
$R$ is reducible to~$S$ iff $\fii\in[\btc+\Gamma,S]$ for all instances $\Gamma\ru\fii$ of~$R$.
\noproof
\end{enumerate}
\end{Obs}

More generally, we say that $R$ is reducible to~$S$ \emph{over a theory~$B$} if $[T,R]\sset[T,S]$ for all $T\Sset B$,
and similarly for weak reducibility. Observation~\ref{obs:rule-red} also generalizes to this situation.

We remark that just like sets of axioms are represented uniquely up to equivalence by \emph{theories}, rules can be
represented up to weak equivalence by \emph{finitary consequence relations,} extending the standard first-order
consequence relation of~$\btc$.

Aside from bounded arithmetic, we will also assume (especially in Section~\ref{sec:pps}) familiarity
with basic propositional proof complexity, and in particular with the quantified propositional sequent calculus~$\G$
(see \cite{book,cook-ngu}). The classes $\siq i$ and~$\piq i$ of quantified propositional formulas are defined as
usual: $\siq0=\piq0$ consists of quantifier-free formulas; $\siq{i+1}$ and~$\piq{i+1}$ include $\siq i\cup\piq i$, and
are closed under $\land$ and~$\lor$; $\siq{i+1}$ is closed under existential quantifiers, and $\piq{i+1}$ under
universal quantifiers; negations of $\siq{i+1}$~formulas are $\piq{i+1}$, and vice versa.

Following~\cite{cook-ngu}, we define $\G_i$ for~$i>0$ as $\G$ restricted so that all cut-formulas are $\siq i$. When
the sequent to be proved consists of $\siq i$~formulas, this is equivalent to the original definition as
in~\cite{book}. Note that up to polynomial simulation, we could allow $\piq i$ cut-formulas in~$\G_i$ as well; on the
other hand, we could restrict cut-formulas to \emph{prenex} $\siq i$~formulas only~\cite{ej-ngu:strict}. Let $\G_i^*$
denote the tree-like version of~$\G_i$. For $i=0$, we define $\G_0$ as extended Frege, optionally considered as a proof
system for prenex $\siq1$ formulas (the system introduced as~$\M{ePK}$ in~\cite{cook-ngu}).

If $P$ is a quantified propositional proof system, and $j\ge0$, then $\rfn_j(P)$ denotes the $\siq j$-reflection
principle for~$P$. If $j=0$, we take this to mean the $\pih1$~reading of the principle: ``for every proof of a
quantifier-free formula~$A$, and every evaluation of subformulas of~$A$ that respects the connectives, the value assigned
to~$A$ is~$1$'' ($\piq0\text-\rfn_P$ in the notation of~\cite[\S X.2.3]{cook-ngu}). (This can make a difference, as
$\btc$ does not necessarily prove that any given quantifier-free formula can be evaluated.) Note that for all proof
systems we are going to consider, this form of $\rfn_0$ is $\btc$-provably equivalent to consistency.

\section{Main systems}\label{sec:main-fragments}

We are ready to introduce the main axioms and rules that will be the topic of this paper. In the rest of this section,
we will show their basic properties, most importantly reductions (inclusions) among the rules.
\pagebreak[2]
\begin{Def}\label{def:main}
The \emph{induction} and \emph{polynomial induction} axioms for a formula~$\fii$ are defined as usual:
\begin{align}
\tag{$\fii$-\ind} \fii(0,y)\land\forall x\,(\fii(x,y)\to\fii(x+1,y))\to\forall x\,\fii(x,y),\\
\tag{$\fii$-\pind} \fii(0,y)\land\forall x\,(\fii(\half x,y)\to\fii(x,y))\to\forall x\,\fii(x,y).
\end{align}
The corresponding induction \emph{rules} are
\begin{gather*}
\tag{$\fii$-\indr} \frac{\fii(0,y)\quad\fii(x,y)\to\fii(x+1,y)}{\fii(x,y)},\\
\tag{$\fii$-\pindr} \frac{\fii(0,y)\quad\fii(\half x,y)\to\fii(x,y)}{\fii(x,y)}.
\end{gather*}
If $\Gamma$ is a set of formulas (usually $\Gamma=\sih i$ or $\Gamma=\pih i$ for $i\ge0$), we define the schema
$\Gamma\text-\ind=\{\fii\text-\ind:\fii\in\Gamma\}$, and similarly for \pind\ and \ppindr.

In the formulation above, the variable~$y$ is a \emph{parameter} of these axioms and rules (we could equivalently allow
a tuple of parameters, as this can be encoded by a single parameter using a pairing function). The corresponding
\emph{parameter-free} schemes, denoted by superscript~${}^-$, are obtained by disallowing~$y$, i.e., $\Gamma$-\ppindf\
consists of $\fii$-\ppind\ for formulas $\fii\in\Gamma$ with no free variables besides~$x$.

The familiar theories $S^i_2$ and~$T^i_2$ are defined as $\btc+\sih i\text-\pind$ and $\btc+\sih i\text-\ind$,
respectively.
\end{Def}

\begin{Rem}
The cases $i=0$ of our schemes and rules are idiosyncratic in various ways: first, $\sih0=\pih0$; second, $\sih0$ is
closed under neither bounded existential nor bounded universal quantifiers, which is going to break some constructions;
and third, $\sih0$-\pind\ and their parameter-free and rule variants are already derivable in the base theory~$\btc$ (that is,
in our language, $S^0_2=\btc$, whereas $T^0_2$ is essentially~$\tpv$).

The standard theories with parameters $T^i_2$ and~$S^i_2$ are axiomatizable by bounded formulas (i.e.,
$\forall\sig\infty$~sentences), since the \ind\ axiom as stated above is equivalent to
\[\forall z\,\bigl(\fii(0,y)\land\forall x<z\,(\fii(x,y)\to\fii(x+1,y))\to\fii(z,y)\bigr),\]
and similarly for~\pind. The proof of this equivalence uses $z$ as a parameter, hence it is not obvious that this
should hold for the parameter-free schemes as well. Nevertheless, the $\pih i$-\ppindf\ schemes do have, for $i>0$,
bounded axiomatizations (specifically, by $\forall\sih{i+1}$~sentences), similarly to the case with parameters: if
$\fii\in\pih i$, then
\begin{equation}\label{eq:15}
\forall x\,\bigl(\fii(0)\land\forall y<x\,(\fii(y)\to\fii(y+1))\to\fii(x)\bigr)
\end{equation}
is provable by induction on the $\pih i$ formula $\psi(x)=\forall y\le x\,\fii(y)$, as
\[\vdash\forall y<x\,(\fii(y)\to\fii(y+1))\land\neg\fii(x)\to\forall z\,(\psi(z)\to\psi(z+1)),\]
and similarly for \pind. This argument does not seem to work for $\sih i$-\ppindf, though.
\end{Rem}

A crucial property is that induction \emph{rules} are equivalent to their parameter-free versions. The case of $\sih
i$ was already proved in~\cite{cffmlm}, but we include it for completeness anyway.
\pagebreak[2]
\begin{Lem}\label{lem:parfreerules}
If\/ $\Gamma=\sih i$ or~$\pih i$ for~$i\ge0$, then $\Gamma\text-\ppindr\equiv\Gamma\text-\ppind^{R-}$.
\end{Lem}
\begin{Pf}
Let $\p{x,y}$ be a $\tc$~pairing function nondecreasing in~$x$ such that $\p{x,y}\ge x+y$, provably in~$\btc$. If $i\le
j\le\lh x$, let $x_{[i,j)}$ denote the number whose binary representation consists of the $i$th through 
$(j-1)$th binary digits of~$x$, where the \emph{most significant} digit has index~$0$; i.e., $x_{[i,j)}=\tdive x{\lh
x-j}\bmod2^{j-i}$.

An instance of $\sih i$-\indr\ for a formula~$\fii(x,y)$ can be reduced to $\sih i$-\indr\ for the formula
$z=0\lor\fii(z_{[m,\lh z)},z_{[1,m)})$, where $m=\cl{\lh z/2}$: we have either $z_{[m,\lh z)}=0$, or $\lh
z=\lh{z-1}$, $z_{[m,\lh z)}=(z-1)_{[m,\lh z)}+1$, and $z_{[1,m)}=(z-1)_{[1,m)}$.

Since $\pih0=\sih0$ and $\btc\vdash\sih0\text-\pind$, we may assume $i>0$ in the remaining cases.

For $\pih i$-\indr, let $\fii(x,y)\in\pih i$, and put $\psi(z)=\forall x,y\le z\,(\p{x,y}\le z\to\fii(x,y))$. Then
\begin{align*}
\fii(0,y)&\vdash\psi(0),\\
\fii(0,y),\fii(x,y)\to\fii(x+1,y)&\vdash\psi(z)\to\psi(z+1),\\
&\vdash\psi(\p{x,y})\to\fii(x,y).
\end{align*}
For $\pih i$-\pindr, we may use $\psi(z)=\forall u\le\lh z\,\fii(z\bmod2^u,\tdive zu)$ in a similar fashion. In order
to verify
\[\fii(0,y),\fii(\half x,y)\to\fii(x,y)\vdash\psi(\half z)\to\psi(z),\]
assume $z>0$, and let $u\le\lh z$. Put $x=z\bmod2^u$, $y=\tdive zu$. If $u=0$, we have $x=0$, and $\fii(0,y)$ holds by
assumption. Otherwise put $z'=\half z$, $u'=u-1$, $x'=z'\bmod2^{u'}$, and $y'=\tdive{z'}{u'}$. We have $u'\le\lh{z'}$,
$x'=\half x$, and $y'=y$, hence $\fii(\half x,y)$ by the induction hypothesis, which implies $\fii(x,y)$ by assumption.

For $\sih i$-\pindr, let $\fii(x,y)$ be a $\sih i$ formula of the form $\exists u\le t(x,y)\,\tet(x,y,u)$ with
$\tet\in\pih{i-1}$. Fix a suitable sequence encoding with $(w)_i$ being the $i$th element of the sequence coded
by~$w$, and $b(z)$ a term such that every sequence~$w$ of length at most~$\lh z$, each of whose entries is bounded by
$t(x,y)$ for some $x,y\le z$, satisfies $w\le b(z)$. Let $\psi(z)$ be the $\sih i$ formula
\[\exists w\le b(z)\,\forall i,j\le\lh z\,\bigl(\p{i,j}<\lh z\to
  (w)_{\p{i,j}}\le t(z_{[j,i+j)},z_{[0,j)})\land\tet(z_{[j,i+j)},z_{[0,j)},(w)_{\p{i,j}})\bigr).\]
Again, the least obvious property to check is that assuming the premises of $\sih i$-\pindr\ for~$\fii$, we can derive
$\psi(\half z)\to\psi(z)$. Let $z>0$, $z'=\half z$, and assume that $w'$ is a sequence of length~$\lh{z'}$
witnessing~$\psi(z')$. We will construct a sequence~$w$ witnessing~$\psi(z)$. If $\p{i,j}<\lh{z'}=\lh z-1$, then
$i+j<\lh{z'}$, thus $z'_{[j,i+j)}=z_{[j,i+j)}$ and $z'_{[0,j)}=z_{[0,j)}$, and we may take
$(w)_{\p{i,j}}=(w')_{\p{i,j}}$. If $\p{i,j}=\lh{z'}$, put $x=z_{[j,i+j)}$, $y=z_{[0,j)}$. Either $i=0$, in which case
$x=0$ and $\fii(0,y)$ holds, or $\p{i-1,j}<\p{i,j}$, $z'_{[0,j)}=y$, and $z'_{[j,j+i-1)}=\half x$. We have $\fii(\half
x,y)$ as witnessed by~$(w)_{\p{i-1,j}}$, hence~$\fii(x,y)$. Either way, we can extend $w'$ to~$w$ so that $(w)_{\p{i,j}}$
is a witness for~$\fii(x,y)$, and then $w$ witnesses~$\psi(z)$.
\end{Pf}
\begin{Cor}\label{cor:indr-indf}
If\/ $\Gamma=\sih i$ or~$\pih i$ for~$i\ge0$, then $\btc+\Gamma\text-\ppindf$ is the weakest theory all of whose
extensions are closed under $\Gamma$-\ppindr.
\end{Cor}
\begin{Pf}
On the one hand, it is clear that any extension of $\Gamma$-\ppindf\ derives $\Gamma\text-\ppind^{R-}$, hence
$\Gamma$-\ppindr\ by Lemma~\ref{lem:parfreerules}. On the other hand, assume that all extensions of~$T$ are closed under
$\Gamma$-\ppindr. Let $\fii\to\psi$ be any instance of $\Gamma$-\ppindf\ as in Definition~\ref{def:main} (here, $\fii$ and~$\psi$ are sentences). Then
$\fii\ru\psi$ is an instance of $\Gamma$-\ppindr, thus $T+\fii\vdash\psi$ by assumption. The deduction theorem then
gives $T\vdash\fii\to\psi$.
\end{Pf}

The next result presents all reductions between our core rules that we know about; they are summarized in
Fig.~\ref{fig:rules}. We will argue in Section~\ref{sec:separations} that no other reductions are likely waiting to be
discovered.
\begin{Thm}\label{prop:basicred}
Let $i\ge0$, and $\Gamma$ be $\sih i$ or~$\pih i$.
\begin{enumerate}
\item\label{item:3}
$\Gamma\text-\ppindr\le\Gamma\text-\ppindf\le\Gamma\text-\ppind$.
\item\label{item:1}
$\sih i\text-\ppind\equiv\pih i\text-\ppind$.
\item\label{item:4}
$\pih i\text-\ppindf\le\sih i\text-\ppindf$, and $\pih i\text-\ppindr\le\sih i\text-\ppindr$.
\item\label{item:5}
$\Gamma\text-\pind\le\Gamma\text-\ind$, $\Gamma\text-\pindf\le\Gamma\text-\indf$, and
$\Gamma\text-\pindr\le\Gamma\text-\indr$.
\item\label{item:6}
$\sih i\text-\ind\le\sih{i+1}\text-\pindr$.
(See also Corollary~\ref{cor:consrules}.)
\item\label{item:7}
$\sih i\text-\indf\le\pih{i+1}\text-\pindf$, and $\sih i\text-\indr\le\pih{i+1}\text-\pindr$.
\end{enumerate}
\end{Thm}
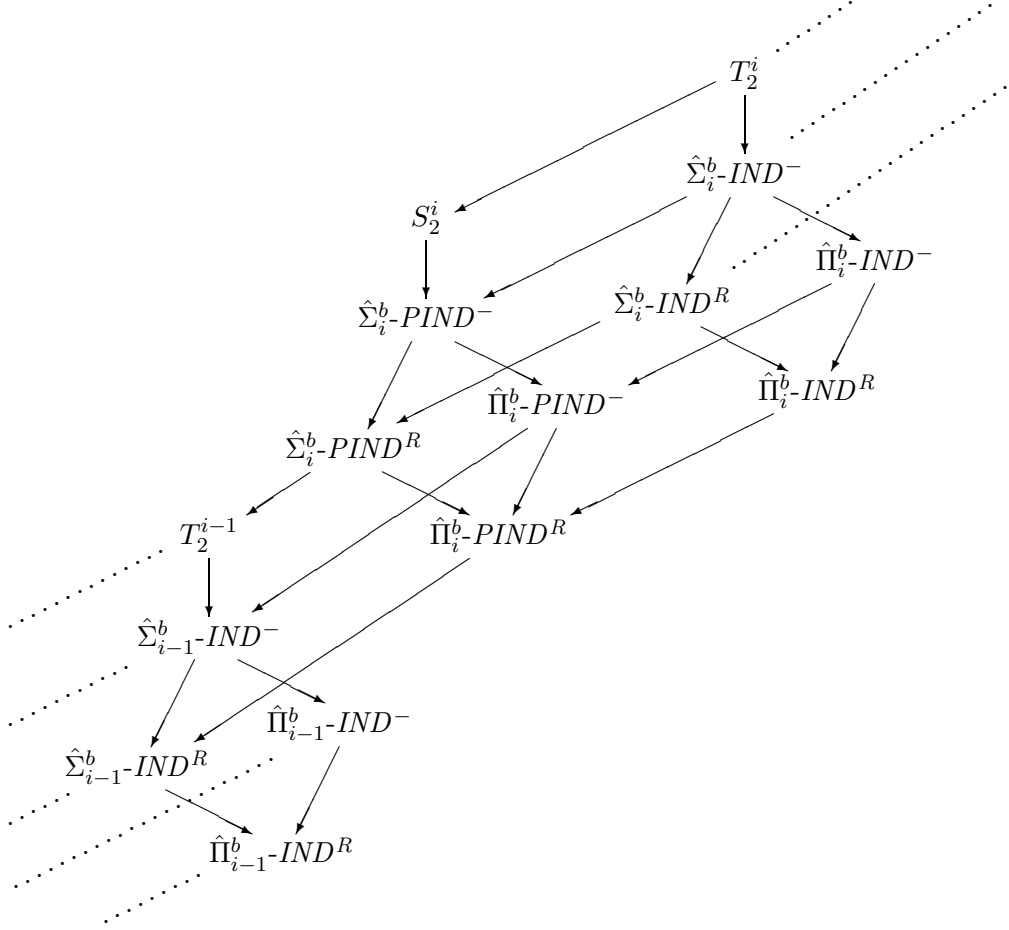
\begin{figure}[tb]
\centering
\magicparoff
\begin{picture}(35,32)(-32,-11)
\multiput(-5.5,19)(.36,.24){8}{.}
\multiput(-5,16)(.36,.24){20}{.}
\multiput(-7,11.5)(.36,.24){27}{.}

\cput(-6.5,18){$T^i_2$}
\put(-6.5,17.5){\vector(0,-1){2}}
\cput(-6.5,14.5){$\sih i$-\indf}
\put(-5.5,14){\vector(2,-1){3}}
\cput(-2,11.5){$\pih i$-\indf}
\put(-7,14){\vector(-1,-2){1.5}}
\cput(-9,10){$\sih i$-\indr}
\put(-2,11){\vector(-1,-2){1.5}}
\put(-8,9.5){\vector(2,-1){3}}
\cput(-4,7){$\pih i$-\indr}

\put(-7.5,18){\vector(-2,-1){9}}
\put(-8.5,14){\vector(-2,-1){7}}
\put(-3.5,11){\vector(-2,-1){7}}
\put(-11.5,9.7){\vector(-2,-1){7}}
\put(-5.5,6.5){\vector(-2,-1){7}}

\cput(-17.5,13){$S^i_2$}
\put(-17.5,12.5){\vector(0,-1){2}}
\cput(-17.5,9.5){$\sih i$-\pindf}
\put(-16.5,9){\vector(2,-1){3}}
\cput(-13,6.5){$\pih i$-\pindf}
\put(-18,9){\vector(-1,-2){1.5}}
\cput(-20,5){$\sih i$-\pindr}
\put(-13,6){\vector(-1,-2){1.5}}
\put(-19,4.5){\vector(2,-1){3}}
\cput(-15,2){$\pih i$-\pindr}

\put(-21.5,4.5){\vector(-3,-2){2.2}}
\put(-14,6){\vector(-3,-2){9.5}}
\put(-16,1.5){\vector(-3,-2){9.5}}

\cput(-25,2){$T^{i-1}_2$}
\put(-25,1.5){\vector(0,-1){2}}
\cput(-25,-1.5){$\sih{i-1}$-\indf}
\put(-24,-2){\vector(2,-1){3}}
\cput(-20.5,-4.5){$\pih{i-1}$-\indf}
\put(-25.5,-2){\vector(-1,-2){1.5}}
\cput(-27.5,-6){$\sih{i-1}$-\indr}
\put(-20.5,-5){\vector(-1,-2){1.5}}
\put(-26.5,-6.5){\vector(2,-1){3}}
\cput(-22.5,-9){$\pih{i-1}$-\indr}

\multiput(-26.8,1.7)(-.4,-.2){14}{.}
\multiput(-28,-2.3)(-.4,-.2){11}{.}
\multiput(-23,-5.5)(-.4,-.2){23}{.}
\multiput(-30,-6.7)(-.4,-.2){6}{.}
\multiput(-25.5,-9.5)(-.4,-.2){9}{.}
\end{picture}
\caption{Reductions between the rules}
\label{fig:rules}
\end{figure}
\begin{Pf}
\ref{item:3} is an immediate consequence of~Lemma~\ref{lem:parfreerules}.

\ref{item:1} is well known: $\ind$ for~$\fii(x,y)$ follows from $\ind$ for~$\neg\fii(a\dotminus x,y)$, and \pind\
for $\fii$ follows from \pind\ for $\neg\fii(\tdive a{\lh x},y)$, where $a$ is an additional parameter.

\ref{item:4}: We may assume $i>0$. Consider an instance of $\pih i$-\indr\ for a formula $\fii(x,y)=\forall z\le t(x,y)\,\tet(x,y,z)$, where
$\tet\in\sih{i-1}$, and let $\psi(x,y,a,z)$ be the $\sih i$ formula 
\[\fii(a\dotminus x,y)\land z\le t(a,y)\to\tet(a,y,z).\]
Then 
\begin{align*}
&\vdash\psi(0,y,a,z),\\
\fii(x,y)\to\fii(x+1,y)&\vdash\psi(x,y,a,z)\to\psi(x+1,y,a,z),\\
\psi(x,y,x,z)&\vdash\fii(0,y)\to\fii(x,y),
\end{align*}
showing that $\fii$-\indr\ reduces to $\psi$-\indr.

In order to show $\pih i\text-\indf\le\sih i\text-\indf$, assume further that $\fii(x)$ is parameter-free. Then
$\btc+\sih i\text-\indf+\fii(0)\land\forall x\,(\fii(x)\to\fii(x+1))$ proves $\fii(x)$ as it is closed under $\sih
i\text-\indr$ by~\ref{item:3}, hence under $\pih i\text-\indr$ by the first part of the proof. Thus, $\btc+\sih
i\text-\indf$ proves $\fii$-\indf\ by the deduction theorem.

The cases of \pindr\ and~\pindf\ are similar, using $\tdive a{\lh x}$ in place of~$a\dotminus x$, as in~\ref{item:1}.

\ref{item:5}: We may assume $i>0$, as $\btc\vdash\sih0\text-\pind$. \pind\ for a $\pih i$~formula~$\fii(x,y)$ follows
from \ind\ for the $\pih i$~formula $\forall u\le x\,\fii(u,y)$, and likewise for \pindf\ or~\pindr. \pind\ for a $\sih
i$~formula $\fii(x,y)$ follows from \ind\ for the formula $\fii(\tdive a{\lh a\dotminus x},y)$ with an additional
parameter~$a$, and this also applies to \pindr. The result for \pindf\ follows from the result for \pindr\ as in the
proof of~\ref{item:4}.

\ref{item:6}: Let
$\fii(x,y)\in\sih i$, and let $\psi(x,y,a)$ be the $\sih{i+1}$ formula
\[\fii(0,y)\land\neg\fii(a,y)\to\exists u\le a,v\le\cl{a/2^{\lh x}}\,(u+v\le a\land\fii(u,y)\land\neg\fii(u+v,y)).\]
Then it is easy to check that $\btc$ proves
\begin{gather*}
\psi(0,y,a),\\
\psi(\half x,y,a)\to\psi(x,y,a),\\
\psi(a,y,a)\to\bigl(\fii(0,y)\land\forall u<a\,(\fii(u,y)\to\fii(u+1,y))\to\fii(a,y)\bigr),
\end{gather*}
thus $[\btc,\sih{i+1}\text-\pindr]$ derives the induction axiom for~$\fii$.

\ref{item:7}: Let $\fii(x,y)\in\sih i$, and let $\psi(x,y,z)$ be the $\pih{i+1}$ formula
\[\forall x'\le z\,(\fii(x',y)\land x+x'\le z\to\fii(x+x',y)).\]
Then
\begin{align*}
&\vdash\psi(0,y,z),\\
\fii(x,y)\to\fii(x+1,y)&\vdash\psi(1,y,z),\\
&\vdash\psi(x_0,y,z)\land\psi(x_1,y,z)\to\psi(x_0+x_1,y,z),\\
\psi(x,y,x)&\vdash\fii(0,y)\to\fii(x,y),
\end{align*}
whence $\sih i\text-\indr\le\pih{i+1}\text-\pindr$. The result for $\sih i$-\indf\ follows as in~\ref{item:4}.
\end{Pf}

\begin{Rem}\label{rem:instances}
Recall that we defined $[T,R]_n$ by counting the nesting depth of applications of~$R$, which is in general necessary in
order to make $[T,R]_n$ a deductively closed first-order theory. However, observe that unnested applications of
\ppindr\ for formulas $\fii_0(x,\vec y)$, \dots, $\fii_k(x,\vec y)$ may be reduced to a single application of the same
rule for the formula $\fii(x,\vec y)=\ET_{i\le k}\fii_i(x,\vec y)$. It follows that if $\Gamma$ is closed
under~$\land$ (such as $\sih i$ or~$\pih i$), then $[T,\Gamma\text-\ppindr]_n$ coincides with the set of formulas
provable using $n$~\emph{instances} of $\Gamma$-\ppindr; the same applies to $\ppind^{R-}$.

Surprisingly, a simple argument shows that the closure of~$T$ under $\pih i$-\ppindr\ collapses to unnested
applications of the rule (thus a single application is enough to prove any given consequence) under very mild
assumptions on the complexity of the theory~$T$. In particular, note that all traditional subsystems of~$S_2$ such as
$S^i_2$ are axiomatized by $\forall\sig\infty\sset\piip$ sentences.
\end{Rem}
\begin{Thm}\label{prop:singlepihindr}
If $T$ is $\piip$-axiomatized, and $i>0$, then
\[T+\pih i\text-\ppindr=[T,\pih i\text-\ppindr].\]
\end{Thm}
\begin{Pf}
In view of Remark~\ref{rem:instances}, it is enough to show that $[T,\pih i\text-\ppindr]$ includes all formulas provable
using two instances of $\pih i\text-\ppind^{R-}$: this implies $[T,\pih i\text-\ppindr]=[T,\pih i\text-\ppindr]_2$,
i.e., $[T,\pih i\text-\ppindr]$ is closed under $\pih i$-\ppindr, and as such it equals $T+\pih i\text-\ppindr$. So, let $\fii,\psi\in\pih i$ be formulas such that
\begin{align*}
T&\vdash\fii(0),\\
T&\vdash\fii(y)\to\fii(y+1),\\
T+\forall y\,\fii(y)&\vdash\psi(0),\\
T+\forall y\,\fii(y)&\vdash\psi(x)\to\psi(x+1).
\end{align*}
(The case of \pind\ is completely analogous.) Since $\psi(0)$ is a bounded sentence, we may assume it is provable
in~$T$ alone. By Parikh's theorem~\ref{thm:parikh}, there is a constant~$c$ such that
\[T\vdash\forall y\le2^{\lh x^c}\,\fii(y)\to(\psi(x)\to\psi(x+1)).\]
Put 
\[\chi(z)=\forall y\le z\,\fii(y)\land\forall x\le z\,(2^{\lh x^c}+x\le z\to\psi(x)).\]
 Then $T$ proves
$\chi(0)$ and $\chi(z)\to\chi(z+1)$, while $\forall z\,\chi(z)$ implies $\forall x\,\psi(x)$.
\end{Pf}

An analogous result for $\sih i$-\ppindr\ only applies to theories~$T$ of bounded complexity ($\forall\sih i$)---more in
line with our expectations---and it seems to require a considerably more complicated proof, see Theorem~\ref{thm:sihinst}.
However, as noted by the reviewer, a weaker result for $\sih i$-\pind\ with $i>0$ and $T$ restricted to
$\exists\forall\sih{i-1}$ follows from Corollary~\ref{cor:buss} and Theorem \ref{prop:basicred} \ref{item:6}, as $T+\sih
i\text-\pindr=T+T^{i-1}_2=[T,\sih i\text-\pindr]$ (cf.\ Corollary~\ref{cor:consrules}).

\section{Variants}\label{sec:variants}

Induction and polynomial induction axioms in bounded arithmetic have equivalent variants that differ in various
details (see e.g.~\cite[\S5.2]{book}): we may consider the length-induction scheme, variants
of minimization principles, or their dual ``ordinal'' induction axioms, and it is not a priori clear if such variants
are still equivalent without parameters. The corresponding induction rules may be varied even more: e.g., the induction
base case may be moved to the conclusion of the rule (cf.~\cite[\S2]{bekl:indru})).

For completeness, we briefly discuss such variants in this section: fortunately, most of them turn out to be
equivalent to some of the axioms and rules introduced in Section~\ref{sec:main-fragments}, except for a few
pathological cases.

\begin{Def}\label{def:variants}
We consider the following schemes and rules, where $\Gamma$ is a set of formulas, and $\fii$ is taken from $\Gamma$:
\allowdisplaybreaks
\begin{gather*}
\tag{$\Gamma$-\lind} \fii(0,y)\land\forall x\,(\fii(x,y)\to\fii(x+1,y))\to\forall x\,\fii(\lh x,y)\\
\tag{$\Gamma\text-\ind_<$} \forall x\,(\forall x'<x\,\fii(x',y)\to\fii(x,y))\to\forall x\,\fii(x,y)\\
\tag{$\Gamma\text-\lind_<$} \forall x\,(\forall x'<x\,\fii(x',y)\to\fii(x,y))\to\forall x\,\fii(\lh x,y)\\
\tag{$\Gamma\text-\pind_<$} \forall x\,(\forall x'\,(\lh{x'}<\lh x\to\fii(x',y))\to\fii(x,y))\to\forall x\,\fii(x,y)\\
\tag{$\Gamma\text-\pind_{\res}$} \forall x\,(\forall u\le\lh x\,(u>0\to\fii(\tdive xu,y))\to\fii(x,y))\to\forall x\,\fii(x,y)\\
\tag{$\Gamma\text-\Min$} \exists x\,\fii(x,y)\to\exists x\,(\fii(x,y)\land\forall x'<x\,\neg\fii(x',y))\\
\tag{$\Gamma\text-\lhmin$} \exists x\,\fii(x,y)\to\exists x\,(\fii(x,y)\land\forall x'\,(\lh{x'}<\lh x\to\neg\fii(x',y)))\\
\tag{$\Gamma\text-\lind^R$} \fii(0,y),\fii(x,y)\to\fii(x+1,y)\ru\fii(\lh x,y)\\
\tag{$\Gamma\text-\ind_0^R$} \fii(x,y)\to\fii(x+1,y)\ru\fii(0,y)\to\fii(x,y)\\
\tag{$\Gamma\text-\pind_0^R$} \fii(\half x,y)\to\fii(x,y)\ru\fii(0,y)\to\fii(x,y)\\
\tag{$\Gamma\text-\lind^R_0$} \fii(x,y)\to\fii(x+1,y)\ru\fii(0,y)\to\fii(\lh x,y)\\
\tag{$\Gamma\text-\ind_<^R$} \forall x'<x\,\fii(x',y)\to\fii(x,y)\ru\fii(x,y)\\
\tag{$\Gamma\text-\lind_<^R$} \forall x'<x\,\fii(x',y)\to\fii(x,y)\ru\fii(\lh x,y)\\
\tag{$\Gamma\text-\pind_<^R$} \forall x'\,(\lh{x'}<\lh x\to\fii(x',y))\to\fii(x,y)\ru\fii(x,y)\\
\tag{$\Gamma\text-\pind_{\res}^R$} \forall u\le\lh x\,(u>0\to\fii(\tdive xu,y))\to\fii(x,y)\ru\fii(x,y)\\
\tag{$\Gamma\text-\Min^R$} \exists x\,\fii(x,y)\ru\exists x\,(\fii(x,y)\land\forall x'<x\,\neg\fii(x',y))\\
\tag{$\Gamma\text-\lhmin^R$} \exists x\,\fii(x,y)\ru\exists x\,(\fii(x,y)\land\forall x'\,(\lh{x'}<\lh x\to\neg\fii(x',y)))
\end{gather*}
As before, the parameter-free versions of these schemes and rules are denoted by~${}^-$.
\end{Def}

\begin{Prop}\label{prop:variants}
Let $\Gamma=\sih i$ or~$\pih i$, where~$i\ge0$, and $\ob\Gamma$ be its dual. The following equivalences hold:
\begin{align}
\label{eq:8}\Gamma\text-\ppind^{R-}_{(0)}&\equiv\Gamma\text-\ppindr,\\
\label{eq:9}\Gamma\text-\lind^{(R)}_{(<)}&\equiv\Gamma\text-\pind^{(R)},\\
\label{eq:10}\sih i/\pih{i+1}\text-\ppind^{(R)(-)}_<&\equiv\pih{i+1}\text-\ppind^{(R)(-)},\\
\label{eq:11}\Gamma\text-\pind_{\res}^{(R)(-)}&\equiv\Gamma\text-\pind^{(R)(-)},\\
\label{eq:12}\Gamma\text-\sch{(P/L)IND}^R_0&\equiv\sih i\text-\ppindr,\\
\label{eq:13}\Gamma\text-\llmin^{(-)}&\equiv\ob\Gamma\text-\ppind^{(-)}_<,\\
\label{eq:14}\Gamma\text-\llmin^R&\equiv\Gamma\text-\llmin.
\end{align}
\end{Prop}
\pagebreak[2]
\begin{Pf}[(sketch)]

\eqref{eq:8}: The position of $\fii(0)$ is immaterial as it is a bounded sentence, and therefore provable or refutable
in~$\btc$. The rest was proved in Lemma~\ref{lem:parfreerules}.

\eqref{eq:9}: \pind\ for $\fii(x,y)$ can be reduced to \lind\ for $\fii(\tdive z{\lh z-x},y)$, while \lind\ for
$\fii(x,y)$ can be reduced to \pind\ for $\fii(\lh x,y)$. In the case of~$\lind_<$, we may use $\forall u\le\lh
x\,\fii(u,y)$; if $\Gamma=\sih i$ (where w.l.o.g.\ $i>0$), we write $\fii(x,y)=\exists z\le t(x,y)\,\tet(x,y,z)$, and use \pind\ on $\exists w\,\forall
u\le\lh x\,\tet(u,(w)_u,y)$  with a suitable bound on~$w$.

\eqref{eq:10}: $\ppind^{(R)(-)}_<$ for $\fii(x,y)$ follows from $\ppind^{(R)(-)}$ for $\forall z\le x\,\fii(z,y)$. On
the other hand, let $\fii(x)=\forall z<2^{\lh x^c}\,\tet(x,z)$ with $\tet\in\sih i$. Then the pairing function
$\p{u,v}:=u2^{\lh u^c}+v$ satisfies $\p{u,v}<\p{u',v'}$ or $\lh{\p{u,v}}<\lh{\p{u',v'}}$ as long as $u<u'$ or $\lh
u<\lh{u'}$ (resp.), $v<2^{\lh u^c}$, and $v'<2^{\lh{u'}^c}$. Thus, defining $\psi(x)$ as
$r(x)<2^{\lh{l(x)}^c}\to\tet(l(x),r(x))$, where $l(\p{u,v})=u$ and $r(\p{u,v})=v$, $\pih{i+1}\text-\ppind^{(R)-}$
for~$\fii$ reduces to $\sih i\text-\ppind^{(R)-}_<$ for~$\psi$. The case with parameters is similar, but easier.

\eqref{eq:11}: $\pind_{\res}$ for $\fii(x,y)$ reduces to \pind\ for $\forall u\le\lh x\,\fii(\tdive xu,y)$; in the case
of~$\Gamma=\sih i$, we swap the outermost quantifiers as in the proof of \eqref{eq:9}.

\eqref{eq:12}: $\sih i\text-\sch{(P/L)IND}^R_0$ is equivalent to $\pih i\text-\sch{(P/L)IND}^R_0$ as in
Theorem \ref{prop:basicred} \ref{item:1}, and it is provable from $\sih i\text-\sch{(P/L)IND}^R$ by replacing
$\fii(x,y)=\exists z\le t(x,y)\,\tet(x,y,z)$ with $z\le t(0,y)\land\tet(0,y,z)\to\fii(x,y)$. (If $i=0$, we just take $\tet=\fii$.)

\eqref{eq:13}: $\llmin$ for~$\fii(x,y)$ amounts to $\ppind_<$ for $\neg\fii(x,y)$.

\eqref{eq:14}: Since $\sih{i+1}\text-\llmin\equiv\pih i\text-\llmin$
by \eqref{eq:10} and~\eqref{eq:13}, it suffices to show $\pih i\text-\llmin\le\pih
i\text-\llmin^R$. Let $\fii(x,y)\in\pih i$. If $i=0$, put $\tet=\fii$, otherwise write $\fii(x,y)=\forall
v\,\tet(x,y,v)$, where $\tet\in\sih{i-1}$.
Let $\psi(x,y,x_0)$ be the $\pih i$ formula
\[\tet(x_0,y,x)\to\fii(x,y).\]
Then $\btc$ proves $\exists x\,\psi(x,y,x_0)$: either $\neg\tet(x_0,y,x)$ for some~$x$, or $\fii(x_0,y)$ and we may
take $x=x_0$. If $\exists x\,\fii(x,y)$, fix $x_0$ such
that $\fii(x_0,y)$. Then a (length-)minimal $x$ satisfying $\psi(x,y,x_0)$ is a (length-)minimal element
satisfying~$\fii(x,y)$.
\end{Pf}

Proposition~\ref{prop:variants} shows that each schema or rule from Definition~\ref{def:variants} is equivalent to one of those
introduced in Definition~\ref{def:main}, except for the following, which are too weak, and thus do not fit nicely in the main
hierarchy:
\begin{itemize}
\item $\Gamma\text-\lind^{(R)-}_{(0/<)}$: Bounded formulas applied to lengths (without non-length parameters) are essentially sharply
bounded, thus $\lind^-$ (as well as all its variants) for bounded formulas whose bounding terms are polynomials is provable in~$\btc$, and full
$\sig\infty\text-\lind^-$ is provable in $\btc+\Omega_2$.
\item $\Gamma\text-\llmin^{R-}$: The premises and conclusions of these rules are $\Sigma_1$ sentences, hence provable
in~$\btc$ if true. It follows that every $\Sigma_1$-sound theory, and every $\piip$-axiomatized
theory, is closed under these rules.
\end{itemize}
Other common variants of induction axioms include maximization schemes. In the presence of parameters, variants of
maximization are easily seen to be equivalent to the corresponding variants of minimization. However, it is unclear how
to sensibly formulate maximization axioms and rules without parameters: the problem is that unlike minimization, we
need an upper bound for maximization, and if this is given by an extra variable, it can be abused to encode arbitrary
parameters.

\section{Conservation}\label{sec:conserv}

In this section we investigate conservation results between induction schemes with and without parameters and
induction rules. The main results state that for theories $T$ of appropriate complexity, $T+T^i_2$ ($T+S^i_2$) is
conservative over $T+\sih i\text-\ppindr$ and $T+\pih i\text-\ppindr$ w.r.t.\ suitable classes of formulas. This will
also imply certain conservativity of $T^i_2$ ($S^i_2$) over $\sih i$-\ppindf\ and $\pih i$-\ppindf.

We start with the easier, and already understood, case of $\sih i$~rules.
The conservation result for $\sih i$-\ppindr\ below, which also implies a conservation result for $\sih i$-\ppindf, was
proved by Cord\'on-Franco, Fern\'andez-Margarit, and Lara-Mart\'\i n~\cite{cffmlm} by model-theoretic means. It generalizes
the special case for $T\sset\forall\sih i$ shown proof-theoretically by Bloch~\cite{bloch}; an analogous result for
$IE_n^-$ was shown earlier by Kaye~\cite{kaye:dio}. We include a proof-theoretic proof of the result for completeness.
\begin{Thm}[\cite{cffmlm}]\label{thm:conssig}
Let $i\ge0$, and $T$ be $\forall\exists\sih{i+1}$-axiomatized. Then
the theory $T+S^i_2$ is $\forall\sih i$-conservative over $T+\sih i\text-\pindr$, and
$T+T^i_2$ is $\forall\sih i$-conservative over $T+\sih i\text-\indr$.
\end{Thm}
\begin{Pf}
We may formulate $T+S^i_2$ in sequent calculus with quantifier-free initial sequents for axioms of~$\btc$, bounded
quantifier introduction rules, the \pind\ rule 
\begin{equation}\label{eq:1}
\frac{\Gamma,\fii(\half x)\seq\fii(x),\Delta}{\Gamma,\fii(0)\seq\fii(t),\Delta},
\end{equation}
where $\fii\in\sih i$ (possibly with parameters not shown) and $x$ is not free in $\Gamma\cup\Delta$, and for every
axiom of~$T$ of the form $\forall x\,\exists y\,\neg\tet(x,y)$ with $\tet\in\sih i$, the rule
\[\frac{\Gamma\seq\tet(t,y),\Delta}{\Gamma\seq\Delta},\]
where $y$ is not free in $\Gamma$, $\Delta$, or~$t$.
By the free-cut-elimination theorem, every $\sih i$ formula provable in $T+S^i_2$ has a sequent proof which only contains $\sih i$
formulas; in particular, the side formulas $\Gamma\cup\Delta$ in each instance of the \pind\ rule are~$\sih i$. Then we
show by \hbox{(meta-)}induction on the length of the proof that all sequents in the proof (that is, their equivalent formulas)
are provable in $T+\sih i\text-\pindr$. The induction step for~\eqref{eq:1} goes as follows. First, we may replace each
formula $\exists u\le s\,\psi(u)$ in~$\Gamma$ with $v\le s\land\psi(v)$, where $v$ is a fresh variable; this turns all
formulas in~$\Gamma$ into $\pih{i-1}\sset\pih i$ formulas (this transformation is not needed if $i=0$, in which case
$\Gamma\sset\sih0=\pih0$ from the get-go). Thus, we may negate them and move them to the right-hand side. Taking
the disjunction of the side formulas on the right-hand side, we are left with a rule
\[\frac{\fii(\half x)\to\fii(x)\lor\psi}{\fii(0)\to\fii(t)\lor\psi},\]
where $\fii,\psi\in\sih i$, and $x$ is not free in~$\psi$. This follows from an instance of $\sih i\text-\pind^R_0$ for
the formula $\fii(x)\lor\psi$, and it is reducible to $\sih i$-\pindr\ by Proposition \ref{prop:variants}~\eqref{eq:12}.

The argument for~$T^i_2$ is similar.
\end{Pf}
Parikh's theorem gives
\begin{Obs}\label{lem:parikh}
If $T$ is $\piip$-axiomatized, then the $\forall\sih i$- and
$\forall\exists\sih i$-fragments of~$T$ are equivalent, for each~$i>0$.
\noproof\end{Obs}

\begin{Cor}\label{cor:conssigres}
Let $i>0$, and $T$ be $\forall\sih{i+1}$-axiomatized. Then
$T+S^i_2$ is $\forall\exists\sih i$-conservative over $T+\sih i\text-\pindr$, and
$T+T^i_2$ is $\forall\exists\sih i$-conservative over $T+\sih i\text-\indr$.
\noproof\end{Cor}

In order to obtain a similar conservation result for $\pih i$-\ppindr\ (Theorem~\ref{thm:picons}), we will need a different method.
Our starting point is the following witnessing theorem, somewhat reminiscent of the KPT
theorem~\cite{kpt}. In the context of parameter-free schemes, it is related to a conservation result for the
$L\Sigma_n^{-\infty}$ scheme (called $I\Pi_n^{-\infty}$ in Kaye~\cite{kaye:parf}) proved by Kaye, Paris, and
Dimitracopoulos~\cite[Thm~2.2]{kpd}.
\begin{Thm}\label{thm:conservp}
Let $i>0$, $T$ be $\forall\exists\sih i$-axiomatized, and $\fii(x)\in\exists\forall\pih i$. If
$T+T^i_2$ ($T+S^i_2$) proves $\forall x\,\fii(x)$, then there are~$k\in\omega$ and $\pih{i-1}$~formulas
$\tet_1(x_0,x_1),\dots,\tet_k(x_0,\dots,x_k)$ such that
\begin{align}
\label{eq:5}
T&\vdash\fii(x_0)\lor\exists y\,\tet_j(x_0,\dots,x_{j-1},y),\qquad j=1,\dots,k,\\
\label{eq:6}
T&\vdash\ET_{j=1}^k\tet_j(x_0,\dots,x_j)\to\fii(x_0)\lor\LOR_{j,l=1}^k\bigl(x_l\prec x_j\land\tet_j(x_0,\dots,x_{j-1},x_l)\bigr),
\end{align}
where $y\prec x$ denotes $y<x$ ($\lh y<\lh x$, respectively).
\end{Thm}
\begin{Pf}
Let $\{\tet_j:j\ge1\}$ be the list of all $\pih{i-1}$~formulas $\tet(\vec x,y)$ such that
\[T\vdash\fii(x_0)\lor\exists y\,\tet(\vec x,y),\]
enumerated in such a way that the free variables of~$\tet_j$ are among
$x_0,\dots,x_{j-1},y$. Put
\[S=T+\neg\fii(c_0)+\{\tet_j(c_0,\dots,c_j):j\ge1\}+\{c_l\prec c_j\to\neg\tet_j(c_0,\dots,c_{j-1},c_l):j,l\ge1\},\]
where $C=\{c_j:j\in\omega\}$ is a set of fresh constants. If the conclusion of the theorem fails, $S$~is consistent.
Let $U$ be a maximal set of $\forall\sih{i-1}(C)$~sentences consistent with~$S$. Let us
fix a model $M\model S+U$, and put $M_0=\{c_j^M:j\in\omega\}$.
\pagebreak[2]
\begin{Cl}
Let $\tet(x_0,\dots,x_n,y)$ be a $\pih{i-1}$~formula such that $M\model\exists y\,\tet(c_0,\dots,c_n,y)$.
\begin{enumerate}
\item\label{item:38}
There are $m\ge n$ and $\psi\in\forall\sih{i-1}$ such that $M\model\psi(c_0,\dots,c_m)$, and
\[T\vdash\psi(x_0,\dots,x_m)\to\fii(x_0)\lor\exists y\,\tet(x_0,\dots,x_n,y).\]
\item\label{item:39}
There exists~$j$ such that $M\model\tet(c_0,\dots,c_n,c_j)$, and $M\model\neg\tet(c_0,\dots,c_n,c_l)$
for all $l$ such that $c_l\prec c_j$.
\end{enumerate}
\end{Cl}
\begin{Pf*}

\ref{item:38}: If not, then $T+\Th_{\forall\sih{i-1}(C)}(M)+\neg\fii(c_0)+\forall y\,\neg\tet(\vec c,y)$ is consistent.
This theory includes $S+U$, but it also contains the $\forall\sih{i-1}(C)$ sentence $\forall y\,\neg\tet(\vec c,y)$
which is not in~$U$ (being false in~$M$), contradicting the maximality of~$U$.

\ref{item:39}: Write $\psi$ as $\forall y\,\xi(x_0,\dots,x_m,y)$ with $\xi\in\sih{i-1}$, and let $j>m$ be such that
$\tet_j(\vec x,y)$ is equivalent to $\neg\xi(\vec x,y)\lor\tet(\vec x,y)$. Then $M\model\tet_j(c_0,\dots,c_j)$,
which means $M\model\tet(c_0,\dots,c_n,c_j)$ as $M\model\xi(c_0,\dots,c_m,c_j)$. Likewise, $M\model c_l\prec
c_j\to\neg\tet(c_0,\dots,c_n,c_l)$.
\end{Pf*}

By part~\ref{item:39} of the claim, $M_0$ is an $\exists\pih{i-1}$-elementary substructure of~$M$. Since
$S\sset\forall\exists\pih{i-1}(C)$, we obtain $M_0\model S$, in particular $M_0\model T+\neg\forall x\,\fii(x)$.

It remains to show $M_0\model T^i_2$ ($S^i_2$, resp.). If $\tet(\vec c,y)$ is a $\pih{i-1}$~formula with parameters
from~$M_0$ such that $M_0\model\exists y\,\tet(\vec c,y)$, then using $M_0\preceq_{\exists\pih{i-1}}M$ and the claim,
there is $j$ such that $M_0\model\tet(\vec c,c_j)$, and $M_0\model\neg\tet(\vec c,c_l)$ for all $l$ such that $c_l\prec
c_j$. Since all elements of $M_0$ are of the form~$c_l$ for some~$l$, this in fact shows
\[M_0\model\tet(\vec c,c_j)\land\forall y\prec c_j\,\neg\tet(\vec c,y).\]
Thus, $M_0\model\pih{i-1}\text-\llmin$, which is equivalent to $\sih i$-\ppind.
\end{Pf}

As an aside, an analogous argument shows the following property, whose special case with $\fii\in\sih i$ may be
employed to give a yet another alternative proof of Theorem~\ref{thm:conssig}:
\begin{Prop}\label{prop:conservs}
Let $i\ge0$, $T$ be $\forall\exists\sih{i+1}$-axiomatized, and $\fii(x)\in\exists\forall\pih{i+1}$. If
$T+T^i_2$ ($T+S^i_2$) proves $\forall x\,\fii(x)$, then there are~$k\in\omega$ and $\pih i$~formulas
$\tet_1(x_0,x_1),\dots,\tet_k(x_0,\dots,x_k)$ satisfying \eqref{eq:5} and
\[T\vdash\ET_{j=1}^k\tet_j(x_0,\dots,x_j)
    \to\fii(x_0)\lor\LOR_{j=1}^k\bigl(x_j\ne0\land\tet_j(x_0,\dots,x_{j-1},P(x_j))\bigr),\]
where $P(x)$ denotes $x-1$ ($\half x$, respectively).
\end{Prop}
\begin{Pf}
We use the same proof as Theorem~\ref{thm:conservp}, with $i'=i+1$ in place of~$i$, and with axioms
\[c_j=0\lor\neg\tet_j(c_0,\dots,c_{j-1},P(c_j))\]
in place of $c_l\prec c_j\to\neg\tet_j(c_0,\dots,c_{j-1},c_l)$ in~$S$. By the same argument, $M_0$ is an
$\exists\pih{i'-1}$-elementary substructure of~$M$ (in particular, $M_0\model T+\neg\forall x\,\fii(x)$), and
$M_0\model\sih{i'-1}\text-\ppind$.
\end{Pf}

\begin{Rem}\label{rem:multidim-min}
The conclusion of Theorem~\ref{thm:conservp} (and, similarly, Proposition~\ref{prop:conservs}) implies that $T$ proves
\begin{multline}\label{eq:16}
\Bigl[\ET_{j=1}^k\forall x_1,\dots,x_{j-1}\,\exists y\,\tet_j(x_0,\dots,x_{j-1},y)\\
          \to\exists x_1,\dots,x_k\,\ET_{j=1}^k\bigl(\tet_j(x_0,\dots,x_j)\land
               \forall z\prec x_j\,\neg\tet_j(x_0,\dots,x_{j-1},z)\bigr)\Bigr]
\to\fii(x_0),
\end{multline}
which means that $\fii(x_0)$ follows over~$T$ from a form of $k$-times iterated $\pih{i-1}$-minimization.

This $k$-dimensional minimization is, similarly to Kaye's $I\Pi_n^{-(k)}$, a form of induction over the
ordinal~$\omega^k$, in contrast to the usual induction over~$\omega$; this is what makes $I\Pi_n^{-\infty}$ strictly
stronger than $I\Pi_n^-$. However, we will see next that in our main case of interest, the $\exists y$ quantifiers
above can be bounded by a term~$t(x_0)$. In that case, the induction is really over the ordinal
$a^k$ for $a=t(x_0)$, which is finite, and as such should follow from ordinary induction. We will formalize this
intuition below.
\end{Rem}
\begin{Lem}\label{lem:conservp-bd}
Let $i>0$, $T\sset\forall\sih i$, and $\fii(x)\in\exists\pih i$. If
$T+T^i_2$ ($T+S^i_2$) proves $\forall x\,\fii(x)$, then there are~$k\in\omega$, $\pih{i-1}$~formulas
$\tet_1(x_0,x_1),\dots,\tet_k(x_0,\dots,x_k)$, and a term~$t(x_0)$ such that
\begin{align}
\label{eq:17}
&\vdash y\ge t(x_0)\to\tet_j(x_0,\dots,x_{j-1},y),\qquad j=1,\dots,k,\\
\label{eq:18}
T&\vdash\ET_{j=1}^k\tet_j(x_0,\dots,x_j)\to\fii(x_0)\lor\LOR_{j=1}^k\exists z\prec x_j\,\tet_j(x_0,\dots,x_{j-1},z),
\end{align}
where $y\prec x$ denotes $y<x$ ($\lh y<\lh x$, respectively).
\end{Lem}
\begin{Pf}
We modify the proof of Theorem~\ref{thm:conservp} as follows. Let $\{\p{\tet_j,t_j}:j\ge1\}$ be an enumeration of pairs
$\p{\tet,t}$ where $t(x)$ is a term, and $\tet(\vec x,y)$ is a $\pih{i-1}$~formula of the form $y\ge t(x_0)\lor\dots$.
We define
\[S=T+\neg\fii(c_0)+\{c_j\le t_j(c_0)\land\tet_j(c_0,\dots,c_j)\land\forall z\prec c_j\,\neg\tet_j(c_0,\dots,c_{j-1},z):j\ge1\},\]
and $U$, $M$, and $M_0$ as in Theorem~\ref{thm:conservp}. Since $S+U$ is $\piip$-axiomatized, its validity is preserved
downwards to cuts; thus, in view of the axioms $c_j\le t_j(c_0)$, we may assume that every element of~$M$ is bounded by
a term in~$c_0$.

In the proof of the Claim, there exists a term $t$ such that $M\model\exists y\prec t(c_0)\,\tet(c_0,\dots,c_n,y)$,
hence we may assume w.l.o.g.\ that $\tet$ has the form $y\prec t(x_0)\land\dots$. We change the definition of
$\tet_j(\vec x,y)$ to $y\ge t(x_0)\lor\neg\xi(\vec x)\lor\tet(\vec x,y)$, with $t_j=t$. Then $M$ satisfies $\tet_j(\vec
c,c_j)$, and $\forall z\prec c_j\,\neg\tet_j(\vec c,z)$. Either $\tet(c_0,\dots,c_n,c_j)$, in which case we are
done, or $c_j=t(c_0)$. But in the latter case, we have $\exists y\prec c_j\,\tet_j(\vec c,y)$, a contradiction.

The rest of the proof is as in Theorem~\ref{thm:conservp}.
\end{Pf}
\begin{Lem}\label{lem:collapse}
In Lemma~\ref{lem:conservp-bd}, we may take $k=1$. That is, under the assumptions of the lemma, there is a
$\pih{i-1}$~formula $\tet(x,y)$ and a term~$t(x)$ such that
\begin{align*}
&\vdash y\ge t(x)\to\tet(x,y),\\
T&\vdash\tet(x,y)\to\fii(x)\lor\exists z\prec y\,\tet(x,z).
\end{align*}
\end{Lem}
\begin{Pf}
Let us first consider the case of \ind. Let $k$, $t$, and $\tet_1,\dots,\tet_k$ be as in Lemma~\ref{lem:conservp-bd}. We
may assume w.l.o.g.\ that $t(x)=2^{\lh x^c}-1$ for some constant $c\ge1$. A $k$-tuple $\p{x_1,\dots,x_k}$ where
$x_1,\dots,x_k<2^{\lh x^c}$ may be represented by a number $y<2^{k\lh x^c}$ as 
\begin{equation}\label{eq:19}
y=x_12^{(k-1)\lh x^c}+x_22^{(k-2)\lh x^c}+\dots+x_k.
\end{equation}
With this encoding in mind, we define a $\pih{i-1}$~formula $\tet(x,y)$ by
\[y\ge2^{k\lh x}\lor\ET_{j=1}^k
    \tet_j\left(x,\dive y{(k-1)\lh x^c}\bmod2^{\lh x^c},\dots,\dive y{(k-j)\lh x}\bmod2^{\lh x^c}\right).\]
Work in~$T$, and assume for contradiction
\[\tet(x,y)\land\forall z<y\,\neg\tet(x,z)\land\neg\fii(x).\]
Since $\tet(x,2^{k\lh x^c}-1)$ by~\eqref{eq:17}, we must have $y<2^{k\lh x^c}$. Write $x_0=x$, and let
$x_1,\dots,x_k<2^{\lh x^c}$ be as in~\eqref{eq:19}. By~\eqref{eq:18}, we have
\[\neg\tet_j(x_0,\dots,x_j)\lor\exists z<x_j\,\tet_j(x_0,\dots,x_{j-1},z)\]
for some~$j=1,\dots,k$. However, $\neg\tet_j(x_0,\dots,x_j)$ is impossible because of $\tet(x,y)$, thus let us fix
$z_j<x_j$ such that $\tet_j(x_0,\dots,x_{j-1},z_j)$, and put
\[z=x_12^{(k-1)\lh x^c}+\dots+x_{j-1}2^{(k-j+1)\lh x^c}+(z_j+1)2^{(k-j)\lh x^c}-1,\]
which represents the $k$-tuple $\p{x_1,\dots,x_{j-1},z_j,2^{\lh x^c}-1,\dots,2^{\lh x^c}-1}$. We have
$\tet_l(x_0,\dots,x_l)$ for $l<j$ as $\tet(x,y)$, $\tet_j(x_0,\dots,x_{j-1},z_j)$ by the choice of~$z_j$, and
$\tet_l(x_0,\dots,x_{j-1},z_j,2^{\lh x^c}-1,\dots)$ for $l>j$ by~\eqref{eq:17}, hence $\tet(x,z)$ and $z<y$, a
contradiction.

In the case of $\pind$, we proceed similarly, except that we encode $\p{x_1,\dots,x_k}$ by
\[2^{\lh{x_1}\lh x^{(k-1)c}+\lh{x_2}\lh x^{(k-2)c}+\dots+\lh{x_k}+k\lh x^c}+x_12^{(k-1)\lh x^c}+x_22^{(k-2)\lh x^c}+\dots+x_k,\]
and we define $\tet(x,y)$ to hold if $y\ge2^{\lh x^{kc}+k\lh x^c}$, or if $y$ is a valid encoding of
$\p{x_1,\dots,x_k}$ such that
\[\ET_{j=1}^k\tet_j(x,x_1,\dots,x_j).\]
It is easy to see that if $y$ encodes $\p{x_1,\dots,x_k}$, and $z$ encodes
$\p{x_1,\dots,x_{j-1},z_j,\dots,z_k}$ with $\lh{z_j}<\lh{x_j}$, then $\lh z<\lh y$. Using this property, the
same proof as above shows
\[T\vdash\tet(x,y)\to\fii(x)\lor\exists z\,\bigl(\lh z<\lh y\land\tet(x,z)\bigr)\]
as required.
\end{Pf}

\pagebreak[2]
\begin{Thm}\label{thm:picons}
If $i>0$ and $T$ is $\forall\sih i$-axiomatized, $T+S^{i+1}_2$ ($T+S^i_2$) is $\forall\pih i$-conservative over $[T,\pih i\text-\ppindr]$.
\end{Thm}
\begin{Pf}
$T+S^{i+1}_2$ is $\forall\sih{i+1}$-conservative over $T+T^i_2$ by Corollary~\ref{cor:buss}, hence it suffices to deal with
$T^i_2$ in place of~$S^{i+1}_2$.

Assume that $T+T^i_2$ ($T+S^i_2$) proves $\forall x\,\fii(x)$ with $\fii\in\sih{i-1}$, and let $\tet$ and~$t$ be as in
Lemma~\ref{lem:collapse}. Putting $\psi(x,y)=\fii(x)\lor\neg\tet(x,y)$, we have
\[T\vdash\forall z\prec y\,\psi(x,z)\to\psi(x,y),\]
hence an application of $\sih{i-1}\text-\ppind_<^R$, equivalent to $\pih i\text-\ppindr$ by Proposition~\ref{prop:variants},
yields $\psi(x,y)$. Since $\tet(x,y)$ holds for all sufficiently large~$y$, this implies~$\fii(x)$.
\end{Pf}

Using a similar strategy, we also obtain a $\sih i$ version of Theorem~\ref{prop:singlepihindr}:
\begin{Thm}\label{thm:sihinst}
If $i>0$ and $T$ is $\forall\sih i$-axiomatized, $T+S^{i+1}_2$ ($T+S^i_2$) is $\forall\sih i$-conservative over $[T,\sih i\text-\ppindr]$. In
particular, $T+\sih i\text-\ppindr=[T,\sih i\text-\ppindr]$.
\end{Thm}
\begin{Pf}
Assume $T+T^i_2$ ($T+S^i_2$) proves $\forall x\,\fii(x)$ with $\fii\in\sih i$, and
let $\tet$ and~$t$ be as in Lemma~\ref{lem:collapse}. In the case of $\sih i$-\ind, we put
\[\psi(x,w)=\fii(x)\lor\exists y\le t(x)\,\bigl(w+y\le t(x)\land\tet(x,y)\bigr),\]
and observe
\begin{align*}
&\vdash\psi(x,0),\\
T&\vdash\psi(x,w)\to\psi(x,w+1),\\
&\vdash\psi(x,t(x)+1)\to\fii(x),
\end{align*}
thus $[T,\sih i\text-\indr]\vdash\fii(x)$. In the case of $\sih i$-\pind, we use
\[\psi(x,w)=\fii(x)\lor\exists y\le t(x)\,\bigl(\lh w+\lh y\le\lh{t(x)}\land\tet(x,y)\bigr)\]
in a similar way.
\end{Pf}

As we will see in Corollary~\ref{cor:sihinst0}, Theorem~\ref{thm:sihinst} also holds for $i=0$.

\begin{Cor}\label{thm:conspi}
Let $T$ be $\forall\sih i$-axiomatized.
\begin{enumerate}
\item\label{item:20}
$T+S^i_2$ is $\forall\exists\sih i$-conservative over $[T,\sih i\text-\pindr]$ for $i\ge1$,
and $\forall\exists\sih{i-1}$-conservative over $[T,\pih i\text-\pindr]$ for $i\ge2$.
\item\label{item:21}
$T+S^{i+1}_2$ is $\forall\exists\sih{i+1}$-conservative over $T+T^i_2$, $\forall\exists\sih i$-conservative over
$[T,\sih i\text-\indr]$ for $i\ge1$, and $\forall\exists\sih{i-1}$-conservative over $[T,\pih i\text-\indr]$ for
$i\ge2$.
\end{enumerate}
\end{Cor}
\begin{Pf}
By Observation~\ref{lem:parikh}, Corollary~\ref{cor:conssigres}, Theorems \ref{thm:picons} and~\ref{thm:sihinst}, and
Corollary~\ref{cor:buss}.
\end{Pf}
\begin{Rem}\label{rem:cons}
We can further extend Theorems \ref{thm:picons} and~\ref{thm:sihinst} and Corollary~\ref{thm:conspi} using the observation that if $T+S$ is
$\Gamma$-conservative over $[T,R]$ and $\psi\in\Gamma$, then $T+\neg\psi+S$ is $\Gamma$-conservative over
$[T+\neg\psi,R]$ (as long as $\Gamma$ is closed under disjunction).

For example, if $i\ge1$, then $T+S^{i+1}_2$ ($T+S^i_2$) is $\forall\exists\sih i$-conservative over $[T,\sih
i\text-\ppindr]$ whenever $T\sset\forall\sih i\cup\exists\forall\pih i$.

We will not list explicitly all cases as it would get too unwieldy.
\end{Rem}

We can draw a few conclusions from Theorems \ref{thm:conssig} and~\ref{thm:picons}. First, some of our rules collapse over
sufficiently simple base theories; this is analogous to the fact that $T+I\Pi_{n+1}^R=T+I\Sigma_n^R$ for
$T\sset\Pi_{n+1}$ (Beklemishev~\cite{bekl:indru}).
\begin{Cor}\label{cor:consrules}
Let $i\ge0$, and $T$ be a theory.
\begin{enumerate}
\item\label{item:49}
If $T$ is $\forall\sih i\cup\exists\pih i$-axiomatized, then $T+\pih{i+1}\text-\pindr=T+\sih
i\text-\indr$. If $i>0$, this holds also for $T\sset\forall\sih i\cup\exists\forall\pih i$.
\item\label{item:50}
If $T$ is $\exists\forall\sih i$-axiomatized, then $T+\sih{i+1}\text-\pindr=T+T^i_2$.
\end{enumerate}
\end{Cor}
\begin{Pf}

\ref{item:49}:
$T+\pih{i+1}\text-\pindr$ includes $T+\sih i\text-\indr$ by Theorem~\ref{prop:basicred}. On the other hand,
$T+\pih{i+1}\text-\pindr$ is $\forall\sih i$-axiomatized over~$T$, and included in $T+S^{i+1}_2$, hence it is
included in $[T,\sih i\text-\indr]$ by Corollary~\ref{thm:conspi} and Remark~\ref{rem:cons}.

\ref{item:50}:
Likewise, $T+T^i_2\sset T+\sih{i+1}\text-\pindr\sset T+S^{i+1}_2$, and the $\forall\sih{i+1}$-fragment of $T+S^{i+1}_2$
is included in $T+T^i_2$ by Corollary~\ref{cor:buss}.
\end{Pf}
\begin{figure}[tb]
\centering
\magicparoff
\begin{picture}(32,30.5)(-29,-9.5)
\multiput(-6,19.1)(.24,.36){6}{.}
\multiput(-5.5,18.8)(.4,.2){11}{.}
\multiput(-5,16)(.4,.2){20}{.}
\multiput(-7,11.5)(.4,.2){25}{.}

\cput(-6.5,18){$T^i_2$}
\put(-6.5,17.5){\vector(0,-1){2}}
\cput(-6.5,14.5){$\sih i$-\indf}
\put(-5.2,14){\vector(3,-2){3.5}}
\cput(-1.5,10.5){$\pih i$-\indf}
\put(-7,14){\vector(-1,-2){1.5}}
\cput(-9,10){$\sih i$-\indr}
\put(-2,10){\vector(-1,-2){1.3}}
\put(-8,9.5){\vector(3,-2){3}}
\cput(-4,6.2){$\pih i$-\indr}

\put(-7.5,18){\vector(-2,-1){8}}
\put(-8.5,14){\vector(-2,-1){6}}
\put(-3.5,10){\vector(-2,-1){11.7}}
\put(-11,9.5){\vector(-2,-1){9.7}}
\put(-4.8,5.7){\vector(-2,-1){18}}

\cput(-16.5,13.5){$S^i_2$}
\put(-16.5,13){\vector(0,-1){2}}
\cput(-16.5,10){$\sih i$-\pindf}
\put(-16.5,9.5){\vector(0,-1){5.3}}
\cput(-16.5,3){$\pih i$-\pindf}
\put(-17.5,2.5){\vector(-2,-1){2.8}}
\put(-18.5,9.3){\vector(-2,-3){3}}

\cput(-22,3.5){$T^{i-1}_2$}
\put(-22,3){\vector(0,-1){1.8}}
\cput(-22,0){$\sih{i-1}$-\indf}
\put(-21,-0.5){\vector(3,-2){3.5}}
\cput(-17,-4){$\pih{i-1}$-\indf}
\put(-22.5,-0.5){\vector(-1,-2){1.5}}
\cput(-24.5,-4.5){$\sih{i-1}$-\indr}
\put(-17.5,-4.5){\vector(-1,-2){1.2}}
\put(-23.5,-5){\vector(3,-2){3}}
\cput(-19.5,-8.1){$\pih{i-1}$-\indr}

\multiput(-23.8,3.2)(-.4,-.2){14}{.}
\multiput(-25,-0.8)(-.4,-.2){11}{.}
\multiput(-19.5,-4.7)(-.4,-.2){24}{.}
\multiput(-27,-5.2)(-.4,-.2){6}{.}
\multiput(-20,-8.7)(-.4,-.2){6}{.}
\end{picture}
\caption{Inclusions between the theories}
\label{fig:theories}
\end{figure}
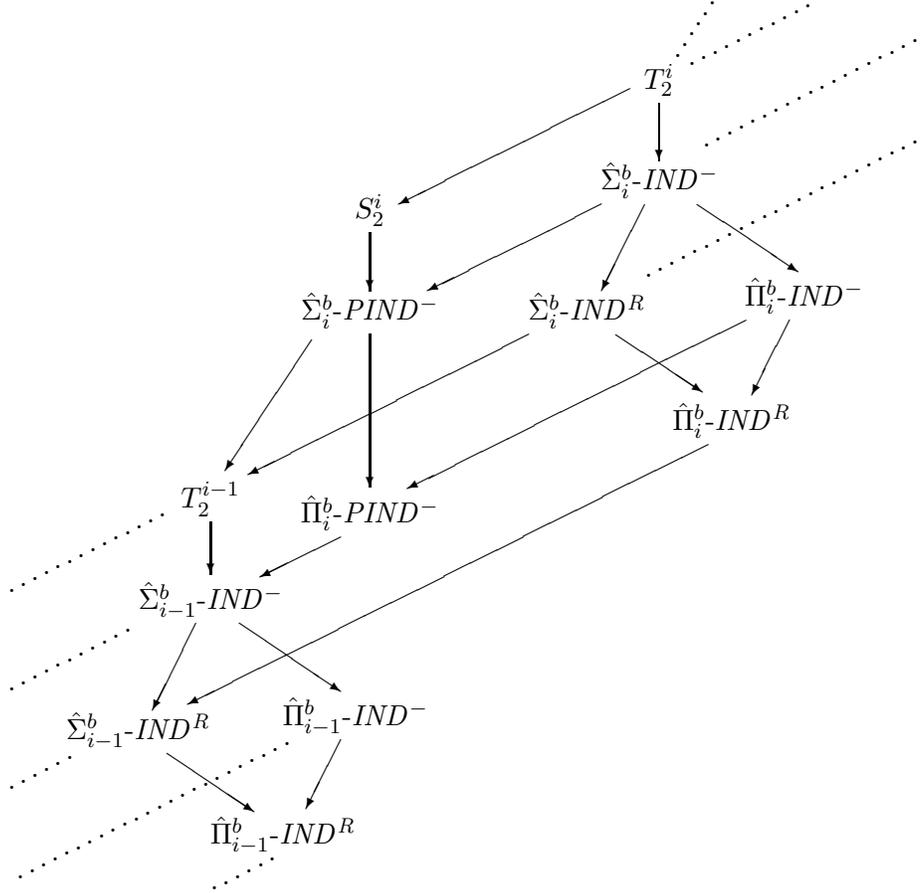

The inclusion diagram between theories axiomatized over~$\btc$ by the rules from Definition~\ref{def:main}, taking into account
Corollary~\ref{cor:consrules}, is depicted in Figure~\ref{fig:theories}. We will present evidence in
Section~\ref{sec:separations} that no further inclusions hold.

Second, we obtain conservation results over parameter-free schemes from the corresponding results for rules and the
deduction theorem. The following corollary summarizes conservativity of $T^i_2$ or~$S^i_2$ over theories axiomatized
over~$\btc$ by parameter-free induction axioms or rules; since the conservations are generally for classes of sentences that
include the complexity of the natural axiomatization of the theories in question, it provides their characterization as
particular fragments of $T^i_2$ or~$S^i_2$.
\begin{Def}\label{def:bool}
If $\Gamma$ is a set of \emph{sentences}, then $\bool(\Gamma)$ denotes the set of Boolean combinations of sentences
from~$\Gamma$, and $\mbool(\Gamma)$ monotone Boolean combinations of sentences from~$\Gamma$.
\end{Def}
\begin{Cor}\label{cor:consindf}
Let $i\ge0$.
\begin{enumerate}
\item\label{item:8} 
$\btc+\sih{i+1}\text-\pindf$ is the $\bool(\forall\sih{i+1})$-fragment of~$S^{i+1}_2$, and it is $\exists\forall\sih{i+1}$-conservative and $\mbool(\exists\pih{i+2}\cup\forall\exists\sih{i+1})$-conservative under~$S^{i+1}_2$.
\item\label{item:2}
$\btc+\sih{i+1}\text-\pindr=T^i_2$ is the $\forall\sih{i+1}$-fragment of~$S^{i+1}_2$, and it is
$\forall\exists\sih{i+1}$-conservative under~$S^{i+1}_2$.
\item\label{item:23}
$\btc+\pih{i+1}\text-\pindf$ is the $\mbool(\exists\pih{i+1}\cup\forall\sih i)$-fragment of~$S^{i+1}_2$, and
if $i>0$, it is $\mbool(\exists\pih{i+1}\cup\forall\exists\sih i)$-conservative under~$S^{i+1}_2$.
\item\label{item:9}
$\btc+\sih i\text-\indf$ is the $\bool(\forall\sih i)$-fragment of $S^{i+1}_2$ or~$T^i_2$, and it is
$\exists\forall\sih i$-conservative under~$T^i_2$. If $i>0$, it is also
$\mbool(\exists\pih{i+1}\cup\forall\exists\sih i)$-conservative under~$T^i_2$, and
$\mbool(\exists\pih i\cup\forall\exists\sih i)$-conservative under~$S^{i+1}_2$.
\item\label{item:10}
$\btc+\sih i\text-\indr=\btc+\pih{i+1}\text-\pindr$ is the $\forall\sih i$-fragment of $S^{i+1}_2$ or~$T^i_2$, and if
$i>0$, it is $\forall\exists\sih i$-conservative under~$S^{i+1}_2$.
\item\label{item:24}
For $i>0$, $\btc+\pih i\text-\indf$ is the $\mbool(\exists\pih i\cup\forall\sih{i-1})$-fragment of $S^{i+1}_2$ or~$T^i_2$.
If $i>1$, it is $\mbool(\exists\pih i\cup\forall\exists\sih{i-1})$-conservative under~$S^{i+1}_2$.
\item\label{item:25}
For $i>0$, $\btc+\pih i\text-\indr$ is the $\forall\sih{i-1}$-fragment of $S^{i+1}_2$ or~$T^i_2$, and if $i>1$, it is
$\forall\exists\sih{i-1}$-conservative under $S^{i+1}_2$.
\end{enumerate}
\end{Cor}
\begin{Pf}
\ref{item:8}:
On the one hand, each instance of $\sih{i+1}$-\pindf\ may be written as an implication between two $\forall\sih{i+1}$
sentences, and it is provable in~$S^{i+1}_2$. On the other hand, if $\fii$ is an $\exists\forall\sih{i+1}$ sentence
provable in~$S^{i+1}_2$, then $\btc+\neg\fii+\sih{i+1}\text-\pindf\Sset\btc+\neg\fii+\sih{i+1}\text-\pindr$ is
inconsistent by Theorem~\ref{thm:conssig} and Lemma~\ref{lem:parfreerules}, thus $\btc+\sih{i+1}\text-\pindf$ proves~$\fii$. Likewise, an
$\mbool(\exists\pih{i+2}\cup\forall\exists\sih{i+1})$ sentence may be written as a conjunction of implications
$\fii\to\psi$, where $\fii\in\forall\sih{i+2}$, and $\psi\in\forall\exists\sih{i+1}$. If $S^{i+1}_2\vdash\fii\to\psi$,
then $\btc+\fii+\sih{i+1}\text-\pindr\vdash\psi$ by Corollary~\ref{cor:conssigres}, thus
$\btc+\sih{i+1}\text-\pindf\vdash\fii\to\psi$.

The other items are similar.
\end{Pf}

Notice that the missing cases $i=0$ in \ref{item:24} and~\ref{item:25} are covered by \ref{item:9} and~\ref{item:10}
(respectively), as $\pih0=\sih0$.

We also obtain similar characterizations for rules rather than theories. For simplicity, we will state only the basic
cases involving $\forall\sig\infty$ formulas. If $\Gamma,\Theta$ are classes of sentences, let us say that a rule $R$ is
a $\Gamma\ru\Theta$ rule if all instances of $R$ have premises in~$\Gamma$ and conclusions in~$\Theta$.
\begin{Cor}\label{cor:consindr}
Let $i\ge0$, and $\Theta=\sih i$ or $\pih i$.
\begin{enumerate}
\item\label{item:47}
$\Theta$-\indr\ is a $\forall\sih i\ru\forall\Theta$ rule derivable in~$T^i_2$, and all $\forall\sih i\ru\forall\Theta$
rules derivable in~$S^{i+1}_2$ are reducible to it.
\item\label{item:48}
$\Theta$-\pindr\ is a $\forall\sih i\ru\forall\Theta$ rule derivable in~$S^i_2$, and all $\forall\sih i\ru\forall\Theta$
rules derivable in~$S^i_2$ are reducible to it.
\end{enumerate}
\end{Cor}
\begin{Pf}
The first assertions in both \ref{item:47} and~\ref{item:48} are clear. Let $\fii\ru\psi$ be a $\forall\sih
i\ru\forall\Theta$ rule derivable in $S^{i+1}_2$ ($S^i_2$, respectively). If $i>0$, we have
$\psi\in[\btc+\fii,\Theta\text-\ppindr]$ by Theorems \ref{thm:picons} and~\ref{thm:sihinst}; if $i=0$, which is a nontrivial case only
for~\ref{item:47}, we have $\psi\in\btc+\fii+\sih0\text-\indr$ by Theorem~\ref{thm:conssig}, and
$\btc+\fii+\sih0\text-\indr=[\btc+\fii,\sih0\text-\indr]$ by Corollary~\ref{cor:sihinst0}. Thus, $\fii\ru\psi$ is reducible to
$\Theta$-\ppindr\ by Observation~\ref{obs:rule-red}.
\end{Pf}

Our third conclusion is that $\sih i$-induction schemes may be extended to variants of $\st\delt{i+1}$-induction.
\begin{Prop}\label{cor:deltaind}
Let $i\ge0$, and $\fii$ be a $\pih{i+1}$ formula.
\begin{enumerate}
\item\label{item:36}
If $\fii$ is provably equivalent to a $\sih{i+1}$ formula in~$S^{i+1}_2$, then $\sih i$-\indf\ proves
$\fii$-\indf, and $\fii\text-\indr$ weakly reduces to $\sih i\text-\indr$; more precisely, $\fii\text-\indr\le\sih
i\text-\indr$ over the theory $[\btc,\sih i\text-\indr]$.
\item\label{item:37}
If $\fii$ is provably equivalent to a $\sih{i+1}$ formula in~$S^i_2$, then $\sih i$-\pindf\ proves
$\fii$-\pindf, and $\fii\text-\pindr$ weakly reduces to $\sih i\text-\pindr$; more precisely,
$\fii\text-\pindr\le\sih i\text-\pindr$ over $[\btc,\sih i\text-\pindr]=T^{i-1}_2$ (over $\btc$ if $i=0$).
\end{enumerate}
\end{Prop}
\begin{Pf}

\ref{item:36}: Let $\fii'$ be a $\sih{i+1}$ formula that $S^{i+1}_2$ proves equivalent to~$\fii$. Recall that
under the assumptions, $\fii$-\ind\ is provable in~$S^{i+1}_2$: assuming $\forall x<a\,(\fii(x,y)\to\fii(x+1,y))$, we
show $\forall x\,(x+z\le a\land\fii'(x,y)\to\fii(x+z,y))$ by $\pih{i+1}$-\pind\ on~$z$.

By Corollary~\ref{cor:consindr}, the $\forall\sih i\ru\forall\sih i$ rule
\[\frac{\fii(0,y)\quad\fii'(x,y)\to\fii(x+1,y)}{\fii(x,y)}\]
is reducible to $\sih i$-\indr, and likewise the $\forall\sih i$ sentence $\forall x,y\,(\fii'(x,y)\to\fii(x,y))$ is
provable in the theory $[\btc,\sih i\text-\indr]$. Thus, $\fii$-\indr\ is derivable from two instances of
$\sih i$-\indr; if $\fii$ is parameter-free, it follows that $\sih i$-\indf\ proves $\fii$-\ind\ by the deduction
theorem.

\ref{item:37} is analogous, using the fact that $S^i_2$ proves $\st\delt{i+1}$-\pind\ \cite[Cor.~8.2.7]{book}. (For
$i=0$, if $\fii$ is $\st\delt1$ in~$\btc$, it is in fact $\sih0$ in~$\btc$, hence $\btc$ proves $\fii$-\pind.)
\end{Pf}

\begin{Rem}
In contrast to Theorem~\ref{thm:conssig}, it is unclear whether the $\forall\pih i$-conservativity of $T+T^i_2$
($T+S^i_2$) over $T+\pih i\text-\ppindr$ in Theorem~\ref{thm:picons}
carries over to $\forall\exists\sih i$-axiomatized theories~$T$, and whether $T^i_2$
($S^i_2$) is $\exists\forall\pih i$-conservative over $\pih i$-\ppindf. (These two problems are in fact equivalent
as a consequence of Theorem~\ref{thm:consindfr} below.)

Notice that the $\forall\pih i$ consequences of $T+T^i_2$ ($T+S^i_2$) are
axiomatized over~$T$ by the rule ``from~\eqref{eq:16} infer $\forall x\,\fii(x)$'' for $\fii\in\sih{i-1}$, and
likewise, the $\exists\forall\pih i$ consequences of $T^i_2$ ($S^i_2$) are axiomatized by the scheme
\begin{multline}\label{eq:2}
\ET_{j=1}^k\forall x_1,\dots,x_{j-1}\,\exists y\,\tet_j(x_1,\dots,x_{j-1},y)\\
          \to\exists x_1,\dots,x_k\,\ET_{j=1}^k\bigl(\tet_j(x_1,\dots,x_j)\land
               \forall z\prec x_j\,\neg\tet_j(x_1,\dots,x_{j-1},z)\bigr)
\end{multline}
for $k\in\omega$ and $\tet_j\in\pih{i-1}$. Thus, the question becomes whether $\pih i$-\ppindf\ proves~\eqref{eq:2}. For
$k=1$, \eqref{eq:2} is just $\pih{i-1}\text-\llmin^-$, which is equivalent to $\pih i$-\ppindf\ by
Proposition~\ref{prop:variants}, hence another formulation is if the scheme~\eqref{eq:2} collapses to its case $k=1$.
\end{Rem}
\begin{Que}\label{que:conservp-unbd}
Let $i>0$.
\begin{enumerate}
\item Is $T^i_2$ ($S^i_2$) $\exists\forall\pih i$-conservative over $\pih i$-\ppindf?
\item Is $T+T^i_2$ ($T+S^i_2$) $\forall\pih i$-conservative over $T+\pih i\text-\ppindr$ for every $\forall\exists\sih
i$-axiomatized theory~$T$?
\end{enumerate}
\end{Que}

Theorems \ref{thm:conssig} and~\ref{thm:picons} imply certain conservativity of \ppindf\ over \ppindr. As we will see below, we can do
better by a direct argument: the conservation results hold over base theories of arbitrary complexity, and they respect
numbers of instances.

Kaye~\cite{kaye:axqc} gave a simple argument showing the conservativity of $k$~instances of axioms of a particular form
over $k$~instances of the corresponding rule, with $I\Sigma_n^R$ as the main intended application. While he states the
result more restrictively, his proof can be seen to give the following general statement.
\begin{Thm}[Kaye~\cite{kaye:axqc}]\label{thm:kayeinst}
Let $\Gamma$ and $\Delta$ be sets of sentences such that $\Gamma\lor\Delta\sset\Gamma$. Let
$A^-=\{\alpha_j\to\beta_j:j<k\}$ be a set of $k$~sentences satisfying $\alpha_j\in\Delta$, and $A^R$ the set of
corresponding rules $\alpha_j\lor\tau\ru\beta_j\lor\tau$ for $\tau\in\Gamma$. Then for any theory $T$, $T+A^-$ is
$\Gamma$-conservative over $[T,A^R]_k$.
\noproof\end{Thm}

Theorem~\ref{thm:kayeinst} implies a conservation result of $\sih i\text-\ppindf$ over $\sih i\text-\ppindr$ preserving
numbers of instances, but it does not seem applicable to $\pih i\text-\ppindr$, as the latter is not invariant
under addition of $\sih i$ side-formulas. We remedy this defect using a modification of Kaye's argument that works under
somewhat different assumptions, at the expense of employing more complicated rules (essentially, several rules from
$A^R$ working in parallel). The conservation results for $\pih i\text-\ppindr$ we proved earlier then allow us to
simulate these rules.
\begin{Lem}\label{lem:varkaye}
Let $\Gamma$ and $\Delta$ be sets of sentences such that $\Gamma\lor\Delta\sset\Gamma$. Let
$A^-=\{\alpha_j\to\beta_j:j<k\}$ be a set of $k$~sentences satisfying $\beta_j\in\Delta$, and let $A^{R\|}$ denote the
rules
\[\frac{\LOR_{j\in J}\alpha_j\lor\tau}{\LOR_{j\in J}\beta_j\lor\tau},\qquad\tau\in\Gamma,J\sset\{0,\dots,k-1\}.\]
Then for any theory~$T$, $T+A^-$ is $\Gamma$-conservative over $[T,A^{R\|}]_k$.
\end{Lem}
\begin{Pf}
Assume that
\begin{equation}\label{eq:7}
T\vdash\ET_{j<k}(\alpha_j\to\beta_j)\to\fii,
\end{equation}
where $\fii\in\Gamma$. We define the sentences
\begin{align*}
\tau_m&=\fii\lor\LOR_{\substack{J\sset k\\\lh J=m}}\ET_{j\in J}\beta_j,\\
\sigma_m&=\fii\lor\LOR_{\substack{J\sset k\\\lh J=m}}\Bigl(\ET_{j\in J}\beta_j\land
  \LOR_{j\notin J}\alpha_j\Bigr)
\end{align*}
for $m\le k$.  Using~\eqref{eq:7}, we can check easily
\begin{align*}
&\vdash\tau_0,\\
&\vdash\sigma_k\to\fii,\\
T&\vdash\tau_m\to\sigma_m,
\end{align*}
it thus suffices to show $[\sigma_m,A^{R\|}]\vdash\tau_{m+1}$. Now, for every $I\sset k$ with $\lh I=k-m$, we have
\[\sigma_m\vdash\fii\lor\LOR_{j\in I}\beta_j\lor\LOR_{j\in I}\alpha_j\]
where $\fii\lor\LOR_{j\in I}\beta_j\in\Gamma$, hence
\[[\sigma_m,A^{R\|}]\vdash\fii\lor\LOR_{j\in I}\beta_j.\]
Since
\[\vdash\ET_{\substack{I\sset k\\\lh I=k-m}}\Bigl(\fii\lor\LOR_{j\in I}\beta_j\Bigr)\to\tau_{m+1},\]
this gives $[\sigma_m,A^{R\|}]\vdash\tau_{m+1}$.
\end{Pf}
\begin{Thm}\label{thm:consindfr}
Let $i\ge0$, and $\Theta=\sih i$ or~$\pih i$. If $T$ is an arbitrary extension of~$\btc$, then $T+\Theta\text-\ppindf$ is
$\forall\Theta$-conservative over $T+\Theta\text-\ppindr$.

More precisely, all $\forall\Theta$ sentences provable from $T$ and $k$~instances of $\Theta$-\ppindf\ are in
$[T,\Theta\text-\ppindr]_k$.
\end{Thm}
\begin{Pf}
We apply Lemma~\ref{lem:varkaye} with $A^-$ being $k$~instances of $\Theta$-\ppindf, and
$\Gamma=\Delta=\forall\Theta$. The rules $A^{R\|}$ are $\forall\sih i\ru\forall\Theta$,
and they are clearly derivable in $T^i_2$ ($S^i_2$, resp.), hence each instance is reducible to an instance of
$\Theta$-\ppindr\ by Corollary~\ref{cor:consindr}.
\end{Pf}
\begin{Cor}\label{cor:finaxindf}
Let $i\ge0$, and $\Theta=\sih i$ or~$\pih i$. If $\Theta$-\ppindf\ is finitely axiomatizable, there is a constant~$k$
such that $T+\Theta\text-\ppindr=[T,\Theta\text-\ppindr]_k$ for every $T\Sset\btc$.
\noproof\end{Cor}

\begin{Que}\label{que:finax}
Are the theories $\sih i$-\ppindf\ and $\pih i$-\ppindf\ finitely axiomatizable?
\end{Que}

\section{Propositional proof systems}\label{sec:pps}

A fundamental tool for analysis of strong theories of arithmetic, especially in the context of induction rules and
parameter-free schemes, are \emph{reflection principles} for other theories of arithmetic
(Beklemishev~\cite{bekl:indru,bekl:parfree}). This idea does not quite work for bounded arithmetic, which is too weak
to prove even the consistency of the base theory~$Q$. Instead, theories of bounded arithmetic may be studied using
reflection principles for \emph{propositional proof systems} by means of translation of bounded formulas to families of
propositional formulas. Apart from the switch from first-order theories to propositional logic, there will be clear
analogies between the form of our results and the classical case of strong systems.

There are two main families of propositional translations of interest:
\begin{enumerate}
\item\label{item:40}
A translation of bounded formulas to \emph{quantified} propositional formulas, where number variables translate to
sequences of propositional variables representing their bits,
and bounded quantifiers translate to blocks of propositional quantifiers.
\item\label{item:41}
A translation of bounded formulas in a \emph{relativized} language (i.e., with a new predicate $\alpha(x)$) to
bounded-depth propositional formulas, where number variables are set to constants, atomic formulas involving~$\alpha$
translate to propositional variables, and bounded quantifiers translate to large disjunctions and conjunctions.
\end{enumerate}
Translation~\ref{item:40} goes back to Cook~\cite{cook} who introduced it as a translation of the equational
theory~$\pv$ to~$\M{EF}$; the extension to quantified propositional logic is due to Kraj\'\i\v cek and
Pudl\'ak~\cite{krpu}. Under this translation, Buss's theories $T^i_2$ correspond to subsystems of the quantified
propositional calculus~$\G$. See Kraj\'\i\v cek~\cite{book} and Cook and Nguyen~\cite{cook-ngu} for detailed
treatments.

Translation~\ref{item:41} was introduced by Paris and Wilkie~\cite{par-wil:counting} for $\idz(\alpha)$. Under this
translation, relativized Buss's theories $T^i_2(\alpha)$ translate to quasipolynomial-size bounded-depth proofs. See
\cite[\S3]{bkz} for a thorough discussion of variants of the Paris--Wilkie translation\footnote{Their setup includes
modular counting gates, but most of the results work also in the usual setup.}.

The relationship between the two translations depends on the point of view. On the one hand, translation~\ref{item:41}
produces exponentially larger formulas than translation~\ref{item:40}. On the other hand, if we identify Buss's
theories with the two-sorted theories $V^i$ using the $\rsuv$-isomorphism, translation~\ref{item:41} becomes
essentially equivalent to a special case of translation~\ref{item:40} for sharply bounded formulas (this is how it
appears in~\cite{cook-ngu}).

In this paper, we are going to work with translation~\ref{item:40}. For one thing, it is already well known that it
leads to an exact correspondence of various subsystems of~$S_2$ (with parameters) to reflection principles for
subsystems of~$\G$, and the setup works smoothly enough so that it can be generalized to the theories we are interested
in.

Perhaps more importantly, translation~\ref{item:41} inherently needs relativized theories, and this is
problematic in the context of parameter-free induction axioms. On the one hand, oracles are somewhat similar to
parameters in that they provide black-box information shared by all parts of the induction axiom, and as such go
against the idea of disallowing parameters; in some contexts, they may be used to sneak parameters back in. See
Section~\ref{sec:relat-separ} for more discussion. On the other hand, the Paris--Wilkie translation~\ref{item:41}
largely eliminates the distinction between induction axioms with and without parameters, as parameters (like all
variables) are set to constants before the translation. This stands in contrast to translation~\ref{item:40}, in which
parameters explicitly manifest as tuples of propositional variables that appear both in premises and conclusions
of translations of induction axioms, and thus their presence makes a difference.

In light of this discussion, for any formula $\fii(\vec x)\in\sig\infty$, let $\{\str\fii_n:n\in\omega\}$ denote a
sequence of quantified propositional formulas obtained by a \ref{item:40}-style translation of~$\fii$, where each
first-order variable~$x_i$ translates to a vector of $n$~propositional variables in~$\str\fii_n$, representing an
integer ${}<2^n$. We do not want to get into the gory technical details of the translation; we can generally follow the
definition of $\|\fii\|^n_{q(n)}$ (for a suitably chosen bounding polynomial~$q(n)$) from Kraj\'\i\v
cek~\cite[\S9.2]{book}, or up to the $\rsuv$~isomorphism, the definition of $\|\fii(\vec X)\|$ in~\cite[\S
VII.5]{cook-ngu}. In particular:
\begin{itemize}
\item bounded existential (universal) quantifiers translate to polynomial-size blocks of existential (universal, resp.)
propositional quantifiers,
\item sharply bounded existential (universal) quantifiers within $\sih0$~formulas translate to polynomial-size
disjunctions (conjunctions, resp.), and
\item propositional connectives translate to themselves.
\end{itemize}
There is a bit of a problem in the definition of the translation for atomic formulas~$\fii$, which we would like to
turn into $\siq0$ (i.e., quantifier-free) formulas: the translation from~\cite{book} is not suitable as it
translates atomic formulas to $\siq1$~formulas (provably equivalent to $\piq1$~formulas in strong enough proof
systems); the translation from~\cite{cook-ngu} does translate atomic (and $\Sig0$) formulas to $\siq0$~formulas---even
of bounded depth---but it only works in a much less expressive language. It does not apply to our $\tc$~language.

The solution is to construct, in a suitably canonical way depending on the exact definition of~$\btc$, for each atomic
formula~$\fii$ a uniform sequence of $\tc$~circuits that compute
it, and expand them into (log-depth) propositional formulas~$\str\fii_n$ by means of formulas computing majority. Something
similar was done in \cite{ej:vnc} for a theory whose language includes $\nci$~functions. Again, the details do not
matter for us, as long as the translation is sufficiently well-behaved so that it can be operated by our theories and
proof systems. We stress that the weakest proof system in which we will reason with the translations is extended Frege.

In this way, the translations of $\sih i$~formulas are $\siq i$, and translations of $\pih i$~formulas are~$\piq i$,
for any $i\ge0$.

We recall the following characterization \cite[X.2.23--24]{cook-ngu} (cf.~\cite{cook-ngu:err}):
\begin{Thm}\label{thm:rfnsiti}
\ \begin{enumerate}
\item If $i\ge j>0$, the $\forall\sih j$ consequences of~$S^i_2$ are axiomatized by $\btc+\rfn_j(\G_i^*)$. If
additionally $i>j$, they are also axiomatized by $\btc+\rfn_j(\G_{i-1})$.
\item If $i>0$, $S^i_2=\btc+\rfn_{i+1}(\G_i^*)$.
\item If $i\ge0$, $T^i_2=\btc+\rfn_{i+1}(\G_{i+1}^*)$.\qedhere
\end{enumerate}
\end{Thm}

The main result of this section will be a characterization of parameter-free induction axioms and induction rules
analogous to Theorem~\ref{thm:rfnsiti}. It will involve the following proof systems:

\begin{Def}\label{def:giax}
Let $i\ge0$. For any $\xi(x)\in\sih i$, we define the proof system $\G_i+\xi$ as $\G_i$
with additional initial sequents of the form $\seq\str\xi_n(\vec A)$, where $n\in\omega$, and $A_0,\dots,A_{n-1}$ are
quantifier-free formulas; $\G_i^*+\xi$ is its tree-like version. 
\end{Def}

\begin{Prop}\label{lem:pps-sim}
Let $i\ge0$, $\xi\in\sih i$, and $\fii\in\sig\infty$.
\begin{enumerate}
\item\label{item:30}
If $i>0$ and $S^i_2+\forall x\,\xi(x)\vdash\forall x\,\fii(x)$, then $\btc$ proves that the formulas $\str\fii_n$ have
$\tc$-constructible polynomial-size $(\G_i^*+\xi)$-proofs.
\item\label{item:31}
If $i>0$ or $\fii\in\sih1$, and $T^i_2+\forall x\,\xi(x)\vdash\forall x\,\fii(x)$, then $\btc$ proves that the formulas
$\str\fii_n$ have $\tc$-constructible polynomial-size $(\G_i+\xi)$-proofs.
\end{enumerate}
\end{Prop}
\begin{Pf}
For $i>0$, the standard proofs of these results without~$\xi$ as in \cite[VII.5.2, X.1.21]{cook-ngu} proceed as
follows. We formulate $S^i_2$ ($T^i_2$) in a sequent calculus with bounded quantifier introduction rules, and an
appropriate induction rule. By the free-cut-elimination theorem, each bounded consequence of the theory has a proof
that only contains bounded formulas such that all cut-formulas are $\sih i$. Then we translate the proof to
propositional logic line by line, supplying short subderivations for each step. This argument works in our situation
just the same: if we enhance the first-order calculus with substitution instances of $\xi\in\sih i$ as additional
axioms, the free-cut-elimination theorem again makes all cuts $\sih i$, and then the same translation as before
produces a valid $\G_i^{(*)}$ proof except for instances of~$\xi$, which translate to the additional axioms of
$\G_i^{(*)}+\xi$. The case $i=0$ needs a different argument (either direct as in \cite[X.1.23]{cook-ngu}, or by
simulation of $\G_1^*$ \cite[VII.4.16]{cook-ngu}), but again it works in the presence of additional quantifier-free
axioms.
\end{Pf}

\begin{Lem}\label{lem:pps-rfn}
Let $i\ge0$, and $\xi\in\sih i$.
\begin{enumerate}
\item\label{item:32}
$T^i_2+\forall x\,\xi(x)$ proves $\rfn_{\max\{i,1\}}(\G_i+\xi)$.
\item\label{item:33}
If $i>0$, $S^i_2+\forall x\,\xi(x)$ proves $\rfn_{i+1}(\G_i^*+\xi)$.
\item\label{item:34}
If $i=0$, $[\btc+\forall x\,\xi(x),\sih0\text-\indr]$ proves $\rfn_0(\G_0+\xi)$.
\end{enumerate}
\end{Lem}
\begin{Pf}
\ref{item:32}: The implication $\forall
x\,\xi(x)\to\rfn_i(\G_i+\xi)$ is $\forall\exists\sih{i+1}$, hence it is enough to prove it in~$S^{i+1}_2$, which is
straightforward for $i>0$: a $(\G_i+\xi)$-proof of a $\siq i$ formula contains only $\siq i$ formulas, hence we may show
by $\pih{i+1}$-\lind\ on the length of the proof that every sequent in the proof is valid. For $i=0$, we may e.g.\ show
that the given assignment can be extended to satisfy all extension axioms in the proof using $\sih1$-\lind, and then
show that all lines of the proof are true under this assignment by $\st\delt1$-\lind. This shows that the target
$\siq1$ formula has a true witness, and therefore is itself true.

\ref{item:33}: We may get rid of each axiom $\seq\str\xi_n(A_0,\dots,A_{n-1})$ in a $(\G_i^*+\xi)$-proof by adding the
$\siq{i+1}$ sentence $\exists x_0,\dots,x_{n-1}\,\neg\str\xi_n(x_0,\dots,x_{n-1})$ to the succedent of every sequent in
the proof. It follows using Theorem~\ref{thm:rfnsiti} that the original end-sequent or one of the new formulas is
true under any given assignment, however, the latter contradicts $\forall x\,\xi(x)$.

\ref{item:34}: It suffices to prove the consistency of $\G_0+\xi$, i.e., $\M{EF}+\xi$. By introducing extension
variables for all subformulas used in the proof and other standard manipulations, $\btc$ knows that if there is an
$(\M{EF}+\xi)$-proof of~$\bot$, there is one where all formulas have bounded size (in particular, we can evaluate them
on any given assignment in~$\tc$), and the only variables that occur in the proof are extension variables. Let
$\pi(z)$ be a $\sih0$ formula stating that $z$ is a proof of this form. Let
\[q_{m-1}\eq A_{m-1},\ q_{m-2}\eq A_{m-2}(q_{m-1}),\ \dots,\ q_0\eq A_0(q_1,\dots,q_{m-1})\]
be the list of all extension axioms used in~$z$. Writing $u_i$ for the $i$th bit of~$u$, let $\tet(u,z)$ be the formula
\begin{multline*}
\pi(z)\to u<2^m\land\forall i<m\,\bigl[\forall j<m\,\bigl(j>i\to u_j=A_j(u_{j+1},\dots,u_{m-1})\bigr)\land u_i=1\\
  \to A_i(u_{i+1},\dots,u_{m-1})=1\bigr].
\end{multline*}
Notice that assuming $\pi(z)$, we can extract $m$ (which is a length) and~$A_i$ from $z$ by a $\tc$ function, hence we
can write $\tet(u,z)$ as a $\sih0$ formula. Clearly, $\btc$ proves $\tet(0,z)$, and $\pi(z)\to\neg\tet(2^m,z)$, that is,
$\forall u\,\tet(u,z)\to\neg\pi(z)$, which in view of the preceding discussion means that
\[\vdash\forall u,z\,\tet(u,z)\to\rfn_0(\G_0+\xi).\]
It thus suffices to verify
\[\forall x\,\xi(x)\vdash\tet(u,z)\to\tet(u+1,z).\]
Assume for contradiction that $\tet(u,z)\land\neg\tet(v,z)$, where $v=u+1$. We must have $\pi(z)$ and $u<2^m$. Let
$i_0\le m$ be the least index of a $0$-bit of~$u$, so that $u_j=v_j$ for $j>i_0$; $u_{i_0}=0$, $v_{i_0}=1$; and $u_j=1$, $v_j=0$
for $j<i_0$. If $v=2^m$, we can show $A_i(x_{i+1},\dots,x_{m-1})=1$ by reverse induction on~$i<m$ (i.e., $\sih0$-\lind,
available in~$\btc$). If $v<2^m$, let $i<m$ be a witness that $\neg\tet(v,z)$, i.e., 
\[\forall j<m\,\bigl(j>i\to v_j=A_j(v_{j+1},\dots,v_{m-1})\bigr)\land v_i=1\land A_i(v_{i+1},\dots,v_{m-1})=0.\]
Since $v_i=1$, this makes $i\ge i_0$. On the other hand, we cannot have $i>i_0$, as then the same would hold for $u$ in
place of~$v$, contradicting $\tet(u,z)$. Thus, $i=i_0$. This implies
\[\forall j<m\,\bigl(j>i_0\to u_j=A_j(u_{j+1},\dots,u_{m-1})\bigr)\land u_{i_0}=0=A_{i_0}(u_{i_0+1},\dots,u_{m-1}).\]
Using $\tet(u,z)$ and $u_j=1$ for $j<i_0$, we then prove $A_j(x_{j+1},\dots,u_{m-1})=1$ for $j<i_0$ by reverse induction
on~$j$ ($\sih0$-\lind\ again), hence in either case,
\[u_j=A_j(u_{j+1},\dots,u_{m-1})\]
for all $j<m$. In other words, the bits of~$u$ taken as an assignment to the $q_j$ variables satisfy all the extension
axioms. Using $\sih0$-\lind\ once more, we show that the assignment in fact satisfies \emph{all} formulas in the proof:
the induction steps for Frege rules follow from the fact that the rules are sound, and the $\str\xi$ axioms are true
because we assume $\forall x\,\xi(x)$. However, the last formula of the proof, $\bot$, is false, which is a
contradiction.
\end{Pf}

\begin{Thm}\label{thm:ppseq}
Let $i\ge0$.
\begin{enumerate}
\item\label{item:15}
$\sih i$-\indr\ is equivalent to the rule
\[\frac{\xi(x)}{\rfn_i(\G_i+\xi)},\qquad\xi\in\sih i.\]
\item\label{item:14}
$\sih i$-\indf\ is equivalent to the scheme
\[\forall x\,\xi(x)\to\rfn_i(\G_i+\xi),\qquad\xi\in\sih i.\]
\item\label{item:17}
For $i>0$, $\pih i$-\indr\ is equivalent to the rule
\[\frac{\xi(x)}{\rfn_{i-1}(\G_i+\xi)},\qquad\xi\in\sih i.\]
\item\label{item:16}
For $i>0$, $\pih i$-\indf\ is equivalent to the scheme
\[\forall x\,\xi(x)\to\rfn_{i-1}(\G_i+\xi),\qquad\xi\in\sih i.\]
\end{enumerate}
If $i>0$, analogous equivalences hold with \pind\ in place of~\ind, and $\G^*_i$ in place of~$\G_i$.
\end{Thm}
\begin{Pf}
\ref{item:14} and~\ref{item:16} follow from \ref{item:15} and~\ref{item:17} and the deduction theorem.

\ref{item:15}: On the one hand, $\forall x\,\xi(x)\ru\rfn_i(\G_i+\xi)$ is a $\forall\sih i\ru\forall\sih i$ rule derivable
in $T^i_2$ by Lemma~\ref{lem:pps-rfn}, hence it reduces to $\sih i\text-\indr$ by Corollary~\ref{cor:consindr}. (If $i=0$,
Corollary~\ref{cor:consindr} depends on Theorem~\ref{thm:ppseq} through Corollaries \ref{cor:finax} and~\ref{cor:sihinst0}.
However, we get a reduction to $\sih0\text-\indr$ directly from Lemma \ref{lem:pps-rfn} \ref{item:34} and
Observation~\ref{obs:rule-red}.)

On the other hand, let $\forall x\,\xi(x)$ be a $\forall\sih
i$ sentence equivalent to $\fii(0)\land\forall x\,(\fii(x)\to\fii(x+1))$, where $\fii\in\sih i$. Then $T^i_2+\forall
x\,\xi(x)$ proves $\forall x\,\fii(x)$, hence by Proposition~\ref{lem:pps-sim}, the formulas~$\str\fii_n$ have short
$(\G_i+\xi)$-proofs, provably in~$\btc$. Consequently, $\btc+\rfn_i(\G_i+\xi)$ proves that $\str\fii_n$ are
tautologies for every length~$n$, which implies $\forall x\,\fii(x)$ by reasoning in~$\btc$. (Note for $i=0$ or the
$\pih1$ cases that even the $\pih1$-definition of validity of $\str\fii_n$ ensures $\forall x<2^n\,\fii(x)$ for
$\fii\in\sih0$: $\btc$ can construct the evaluation of $\str\fii_n$ and its subformulas under a given assignment using
a $\tc$ function, even though it may not prove that propositional formulas can be evaluated in general.)

\ref{item:17} is similar to~\ref{item:15}, and the arguments for $\pind$ and~$\G^*_i$ are analogous.
\end{Pf}

\begin{Cor}\label{cor:finax}
If $i\ge0$, and $T$ is a finitely $\forall\sih i$-axiomatized extension of~$\btc$, then the theories $T+\sih i\text-\ppindr$ and
$T+\pih i\text-\ppindr$ are finitely axiomatizable.

Specifically, if $T=\btc+\forall x\,\xi(x)$ with $\xi\in\sih i$, then
\begin{align*}
T+\sih i\text-\indr&=\btc+\rfn_i(\G_i+\xi),\\
\intertext{and for $i>0$,}
T+\sih i\text-\pindr&=\btc+\rfn_i(\G^*_i+\xi),\\
T+\pih i\text-\indr&=T+\rfn_{i-1}(\G_i+\xi),\\
T+\pih i\text-\pindr&=T+\rfn_{i-1}(\G^*_i+\xi).
\end{align*}
\end{Cor}
\begin{Pf}
The inclusions~$\Sset$ are special cases of Theorem~\ref{thm:ppseq}. On the other hand, $T+\sih i\text-\indr$ is
$\forall\sih i$-axiomatized, and if $T+\sih i\text-\indr\sset T^i_2+\forall x\,\xi(x)$ proves a $\sih i$ formula
$\fii(x)$, then $\btc+\rfn_i(\G_i+\xi)$ proves $\fii(x)$ by the argument in the proof of Theorem~\ref{thm:ppseq}. The other
cases are similar, except that the arguments work just for $\fii\in\sih{i-1}$ if we have only $\rfn_{i-1}$. This is
fine as $T+\pih i\text-\ppindr$ is $\forall\sih{i-1}$-axiomatized over~$T$.
\end{Pf}

Using this characterization, we can extend Theorem~\ref{thm:sihinst} to the case $i=0$:
\begin{Cor}\label{cor:sihinst0}
If $T\sset\forall\sih0$, $T+\sih0\text-\indr=[T,\sih0\text-\indr]$.
\end{Cor}
\begin{Pf}
W.l.o.g., $T$ is finitely axiomatizable, hence we may write $T=\btc+\forall x\,\xi(x)$ with $\xi\in\sih0$. Then
$T+\sih0\text-\indr=\btc+\rfn_0(\G_0+\xi)\sset[T,\sih0\text-\indr]$ by Corollary~\ref{cor:finax} and Lemma~\ref{lem:pps-rfn}.
\end{Pf}

A direct proof of Corollary~\ref{cor:sihinst0} is also possible, but it is not particularly illuminating.

\begin{Rem}\label{rem:i+1}
We could extend the definition of $\G_i+\xi$ to $\xi\in\sih{i+1}$ as follows: write $\xi(x)=\exists y<2^{\lh
x^c}\,\neg\tet(x,y)$ with $\tet\in\sih i$, and let $\G_i+\xi$ denote $\G_i$ augmented by the rule
\[\frac{\Gamma\seq\Delta,\str\tet_{n,n^c}(A_0,\dots,A_{n-1},x_0,\dots,x_{n^c-1})}{\Gamma\seq\Delta},\]
where $A_j$ are quantifier-free, and $x_j$ are not free in $\Gamma$, $\Delta$, or $A_{j'}$; likewise for $\G_i^*+\xi$.
(This is easily seen to be p-equivalent to the original definition if $\xi\in\sih i$.) Proposition~\ref{lem:pps-sim} continues
to hold in this setting, and the proof of Lemma~\ref{lem:pps-rfn} gives $S^{i+1}_2+\forall
x\,\xi(x)\vdash\rfn_i(\G_i+\xi)$. Since this extension does not seem to yield new insights about parameter-free
induction schemes or rules, we skip the details.
\end{Rem}

\section{Separations}\label{sec:separations}

We have seen in the previous sections many results relating subsystems of bounded arithmetic with and without
parameters, but in order for these results to be useful, it would be nice to know that the systems do not collapse:
what if the parameter-free induction schemes are actually equivalent to the usual schemes with parameters, so that
e.g.\ $T^i_2=\sih i\text-\indf$? This would make the investigation of $\indf$ rather pointless. Likewise, since we
spent so much effort on $\pih i$ schemes and rules, we would like to know that they are genuinely distinct from the
corresponding $\sih i$ rules.

In general, we are interested if there are any reductions between our schemes and
rules that do not follow from Theorem~\ref{prop:basicred} (as depicted in Figure~\ref{fig:rules}), and furthermore if there
are any inclusions between the theories generated by our schemes and rules over the base theory that do not follow from
Theorem~\ref{prop:basicred} and Corollary~\ref{cor:consrules} (as depicted in Figure~\ref{fig:theories}).

Checking all the cases naively would be a gargantuan task: we have 10 schemes and rules at each level of the hierarchy,
and we need to consider reductions spanning three levels: e.g., $S^i_2$ is supposed not to be included in
$\btc+\pih{i+1}\text-\indf$, which is two levels higher up, being $\forall\sih i$-conservative under~$S^{i+2}_2$.
However, we do not actually have to consider all possible pairs, as there is a lot of redundancy: for example, we do
not need to check separately that $\btc+\sih i\text-\indf\nvdash T^i_2$, because $T^i_2\Sset S^i_2$,
$\sih i\text-\indf\sset\pih{i+1}\text-\indf$, and we want to make sure that $\btc+\pih{i+1}\text-\indf\nvdash S^i_2$
anyway. Let us put our job into a more formal setting:
\begin{Def}\label{def:check}
A \emph{basis of non-inequalities} of a poset $\p{P,{\le}}$ is a set $B\sset P^2$ such that
\begin{enumerate}
\item $a\nleq b$ for any $\p{a,b}\in B$, and
\item for each $a,b\in P$ such that $a\nleq b$, there is $\p{a',b'}\in B$ such that $a'\le a$ and $b\le b'$.
\end{enumerate}
A \emph{critical pair} of~$P$ is $\p{a,b}\in P$ such that $a\nleq b$, but $a'\le b$ for all $a'<a$, and $a\le b'$ for
all $b'>b$. Observe that any basis of non-inequalities of~$P$ has to include all critical pairs.

Let $\p{P_R,{\le_R}}$ denote the poset with formal elements representing $\btc$ and the axioms and rules $\sih i$-\ind,
$\sih i$-\indf, $\sih i$-\indr, $\pih{i+1}$-\indf, $\pih{i+1}$-\indr, $\sih{i+1}$-\pind, $\sih{i+1}$-\pindf,
$\sih{i+1}$-\pindr, $\pih{i+1}$-\pindf, and $\pih{i+1}$-\pindr\ for~$i\ge0$, and with $\le_R$ being the transitive
reflexive closure of the relation given by Theorem~\ref{prop:basicred}. ($\btc$ is a least element of~$P_R$.)

Let $\p{P_T,{\le_T}}$ be the quotient of $\p{P_R,{\le_R}}$ identifying $\sih{i+1}$-\pindr\ with~$\sih i$-\ind, and
$\pih{i+1}$-\pindr\ with $\sih i$-\indr, for each $i\ge0$.

Beware that neither $P_R$ nor~$P_T$ is a lattice.
\end{Def}
\begin{Lem}\label{lem:check}
Let $\p{P,{\le}}$ be a poset in which all strictly increasing infinite sequences are upwards cofinal, and all strictly
decreasing infinite sequences are downwards cofinal\footnote{In fact, weaker assumptions suffice: it is enough if
$\mathbb Q$, $\omega\sqcup1$, and $\omega^*\sqcup1$ do not embed in~$P$, where $\sqcup$ denotes disjoint union of
posets.}. Then the set of critical pairs is a basis of non-inequalities of~$P$.
\end{Lem}
\begin{Pf}
The assumptions may be restated such that for each $u\in P$, $<$ is well-founded on $\{x\in P:x\nleq u\}$, and converse
well-founded on $\{x\in P:u\nleq x\}$. Thus, given $a\nleq b$, we can find a minimal $a'\le a$ such that $a'\nleq b$,
and then a maximal $b'\ge b$ such that $a'\nleq b'$. Then $\p{a',b'}$ is a critical pair.
\end{Pf}

The critical pairs of $P_R$ and~$P_T$ can be determined by a somewhat tedious, but straightforward computation, chasing
the diagrams in Figures \ref{fig:rules} and~\ref{fig:theories}. We see that $P_R$ and~$P_T$ have common critical pairs
\begin{align*}
    \langle\sih i\text-\pind,&\pih{i+1}\text-\indf\rangle,& \langle\sih0\text-\ind,&\pih1\text-\indf\rangle,\\
   \langle\pih i\text-\pindf,&\pih{i+1}\text-\indr\rangle,&\langle\sih0\text-\indf,&\pih1\text-\indr\rangle,\\
    \langle\pih i\text-\indr,&\sih i\text-\pind\rangle,   &\langle\sih0\text-\indr,&\btc\rangle\\
\intertext{for $i\ge1$. Moreover, $P_R$ has critical pairs}
\langle\pih{i+1}\text-\pindr,&\sih i\text-\ind\rangle
\end{align*}
for $i\ge0$, but we can disregard these: $\pih{i+1}\text-\pindr\le T^i_2$ implies
$T^i_2\vdash\pih{i+1}\text-\pindf$ using the deduction theorem,
hence also $\btc+\pih{i+2}\text-\indr\vdash\pih{i+1}\text-\pindf$, which is an instance of another critical pair. Thus,
we obtain:
\begin{Prop}\label{prop:critp}
If there is a reduction between the schemes and rules from Definition~\ref{def:main} which does not follow from
Theorem~\ref{prop:basicred}, or an additional inclusion between the first-order theories they generate over~$\btc$ not
warranted by Corollary~\ref{cor:consrules}, it implies one of the following:
\begin{align}
\label{item:22}
S^i_2&\vdash\btc+\pih i\text-\indr\text{ for some }i\ge0,\\
\label{item:26}
\pih{i+1}\text-\indf&\vdash S^i_2\text{ for some }i>0,\\
\tag{\ref{item:26}$'$}\label{item:28}
\pih1\text-\indf&\vdash T^0_2,\\
\label{item:27}
\btc+\pih{i+1}\text-\indr&\vdash\pih i\text-\pindf\text{ for some $i>0$, or}\\
\tag{\ref{item:27}$'$}\label{item:29}
\btc+\pih1\text-\indr&\vdash\sih0\text-\indf.
\end{align}
(Recall that in our setup, $S^0_2=\btc$.)
\noproof\end{Prop}

The remaining goal is to convince ourselves that \eqref{item:22}--\eqref{item:29} are likely false, or at least
suspect. We
are not very picky, and do not attempt to devise sophisticated separation arguments optimized for the particular
theories; rather, we are content with any evidence that we did not overlook something in Theorem~\ref{prop:basicred}. We
will present run-of-the-mill separations of two kinds, as commonly done for systems of bounded arithmetic:
separations conditional on plausible complexity-theoretic assumptions, and unconditional separations of relativized
versions of our theories.

\subsection{Unrelativized separations}\label{sec:unrel-separ}

The state of our knowledge does not allow us to disprove even $\btc=S_2$ unconditionally---this would require a
major breakthrough. We thus cannot disprove \eqref{item:22}--\eqref{item:29} either. What we can do instead is to show
that they imply other statements (from computational and proof complexity) that are more commonly recognized as
implausible.

\begin{Thm}\label{prop:si-piindr}
If $S^i_2\vdash\btc+\pih i\text-\indr$, then $T^i_2$ is $\forall\sih{\max\{i-1,0\}}$-conservative over $S^i_2$ (and
thus over~$T^{i-1}_2$ for $i>0$). Consequently:
\begin{enumerate}
\item\label{item:11} If $i=0$, $\tc$-Frege p-simulates $\M{EF}$.
\item\label{item:12} If $i>0$, $\G_i^*$ and $\G_{i-1}$ p-simulate~$\G_i$ with respect to
$\siq{i-1}$ sequents.
\item\label{item:13} If~$i>1$, the game induction principle $\mathrm{GI}_i$ (Skelley and Thapen~\cite{sk-th:gi})
is reducible to~$\mathrm{GI}_{i-1}$.
\end{enumerate}
\end{Thm}
\begin{Pf}
The conservativity of $T^i_2$ over~$S^i_2$ is a consequence of the characterization of $\btc+\pih i\text-\indr$ from
Corollary \ref{cor:consindf} \ref{item:25} (or \ref{item:10} if $i=0$). Then \ref{item:11} and~\ref{item:12} follow by a standard argument: $T^i_2$, hence $S^i_2$
and~$T^{i-1}_2$ by assumption, proves $\rfn_{i-1}(\G_i)$. Thus, $\btc$ proves that the tautologies
$\str{\rfn_{i-1}(\G_i)}_n$ have $\tc$-constructible proofs in $\G_i^*$ and~$\G_{i-1}$, which in turn implies that these
two proof systems p-simulate $\G_i$-proofs of $\siq{i-1}$ sequents. Similarly, \ref{item:13} follows from the fact that
$\mathrm{GI}_i$ is complete for the class of NP-search problems provably total in~$T^i_2$.
\end{Pf}

Recall that $\fp^{\Sigma^P_i[O(g(n)),\mathrm{wit}]}$ denotes the class of total search problems computable by a
polynomial function that makes $O(g(n))$ queries to a witnessing $\Sigma^P_i$ oracle, meaning that for any positive
answer, the oracle also has to produce a witness to the outermost existential quantifier. For any $i>0$, the $\sih{i+1}$-definable search problems provably
total in~$S^i_2$ comprise exactly $\fp^{\Sigma^P_i[O(\log n),\mathrm{wit}]}$, and the $\sih{i+1}$-definable search
problems provably total in~$T^{i-1}_2$ comprise exactly $\fp^{\Sigma^P_i[O(1),\mathrm{wit}]}$ (see e.g.\
\cite[Thm.~VIII.7.17]{cook-ngu}; the original results are due to Kraj\'\i\v cek, Pudl\'ak, and Takeuti~\cite{kpt} and
Kraj\'\i\v cek~\cite{kra:witlog}).
\begin{Thm}\label{prop:pi+1indf-si}
\ \begin{enumerate}
\item\label{item:35}
If $\pih1\text-\indf\vdash T^0_2$, then $\ptime=\tc$.
\item\label{item:18}
If $\pih{i+1}\text-\indf\vdash S^i_2$ for some $i>0$, then $\fp^{\Sigma^P_i[O(\log
n),\mathrm{wit}]}=\fp^{\Sigma^P_i[O(1),\mathrm{wit}]}$, and $\ph=\bool(\Sigma^P_{i+1})$.
\end{enumerate}
\end{Thm}
\begin{Pf}
First, observe that $\pih{i+1}$-\indf\ follows from the set of all true $\forall\sih i$ sentences: it is axiomatized by
sentences of the form $\fii\to\psi$, where $\fii\in\forall\sig\infty$, and $\psi\in\forall\sih i$. If $\fii$ is false,
$\neg\fii$ (and a fortiori $\fii\to\psi$) is provable in~$\btc$, being a true $\Sigma_1$ sentence. Otherwise, $\psi$
is true, hence included in $\Th_{\forall\sih i}(\stm)$.

\ref{item:35}: Every poly-time function~$f$ has a provably total $\sih1$-definition in~$T^0_2$, hence by assumption, in
$\Th_{\forall\sih0}(\stm)$, i.e., in the set of true universal sentences of~$L_{\tc}$. By Herbrand's theorem (and
closure under definitions by cases), $f$ is definable by an $L_{\tc}$-term, i.e., it is a $\tc$-function. In particular, every
poly-time predicate is computable in~$\tc$.

\ref{item:18}: Every $\fp^{\Sigma^P_i[O(\log n),\mathrm{wit}]}$ search problem has a
$\sih{i+1}$-definition provably total in $S^i_2$, hence by assumption, in $\Th_{\forall\sih i}(\stm)$. We claim that,
just like for~$T^{i-1}_2$, the provably total $\sih{i+1}$-definable search problems of $\Th_{\forall\sih i}(\stm)$ are
in $\fp^{\Sigma^P_i[O(1),\mathrm{wit}]}$: if
\[\forall u\,\psi(u)\vdash\forall x\,\exists y\,\fii(x,y),\]
where $\psi\in\sih i$, $\fii\in\sih{i+1}$, and $\stm\model\forall u\,\psi(u)$, we have
\[T^{i-1}_2\vdash\forall x\,\exists y\,\bigl(\neg\psi(y)\lor\fii(x,y)\bigr)\]
(the $T^{i-1}_2$ is not really doing anything for us here). We may bound the $y$ using Parikh's theorem, and then by
the above-mentioned characterization of $\forall\sih{i+1}$~consequences of~$T^{i-1}_2$, we obtain
\[T^{i-1}_2\vdash\forall x\,\bigl(\neg\psi(f(x))\lor\fii(x,f(x))\bigr)\]
for some search problem $f\in\fp^{\Sigma^P_i[O(1),\mathrm{wit}]}$, $\sih{i+1}$-definable in~$T^{i-1}_2$; but the first disjunct cannot happen in the real
world:
\[\stm\model\forall x\,\fii(x,f(x)).\]

Thus, $\fp^{\Sigma^P_i[O(\log n),\mathrm{wit}]}=\fp^{\Sigma^P_i[O(1),\mathrm{wit}]}$. This implies
$\ptime^{\Sigma^P_i[O(\log n)]}=\ptime^{\Sigma^P_i[O(1)]}=\bool(\Sigma^P_i)$, as predicates (i.e., $\{0,1\}$-valued
functions) in $\fp^{\Sigma^P_i[O(1),\mathrm{wit}]}$ are in $\ptime^{\Sigma^P_i[O(1)]}$ (cf.\ \cite[6.3.4--5]{book}).
This in turn implies the collapse of $\ph$ to
$\bool(\Sigma^P_{i+1})$ by Chang and Kadin~\cite{cha-kad}.
\end{Pf}

\begin{Rem}\label{rem:collapse-ph}
The second point of Theorem~\ref{prop:pi+1indf-si} is a variant of the well-known result that $T^{i-1}_2=S^i_2$ implies the
collapse of $\ph$, originally proved in \cite{kpt}, and subsequently improved
in~\cite{buss:coll,zamb:notes,cookkra,ej:hash}. The current state of the art is that $T^{i-1}_2=S^i_2$ implies
$T^{i-1}_2\vdash\ph=\bool(\Sigma^P_i)$ \cite[Cor.~4.7]{ej:hash}, which is a one whole level deeper collapse than in
Theorem~\ref{prop:pi+1indf-si}.

While we did not attempt to check the details, it is not implausible that these improvements also work in the presence
of additional true $\forall\sih i$ axioms; if correct, this would strengthen the conclusion of
Theorem \ref{prop:pi+1indf-si} \ref{item:18} to $\ph=\bool(\Sigma^P_i)$.
\end{Rem}

\begin{Que}\label{que:unrel}
Can we disprove \eqref{item:27} or~\eqref{item:29} under a credible hypothesis?
\end{Que}

\subsection{Relativized separations}\label{sec:relat-separ}

Rather than relying on unproven hypotheses, we may want to look at unconditional separations of relativized theories.
All theories we work with may be relativized in the standard way: we include a new predicate symbol~$\alpha(x)$ in the
language, and extend all schemes to allow the use of~$\alpha$ along with other atomic formulas, but do not include any
axioms to fix its particular values.

Relativization is commonly employed in bounded arithmetic to obtain separation results, exploiting the
fact that we can unconditionally separate various complexity classes in the relativized setting. The usefulness of this
technique of course hinges on our belief that for the classes in question (e.g., levels of the polynomial hierarchy),
noninclusions between their relativized versions truly reflect properties of the original unrelativized classes.
(Relativized bounded arithmetic is also useful in connection with bounded-depth propositional proof systems, as the
Paris--Wilkie translation only makes sense for relativized theories.)

Relativization of parameter-free schemes may seem somewhat more dubious than in the case of usual theories of bounded
arithmetic, as it goes against the spirit of parameter removal: similar to parameters, the oracle provides access to
additional black-box information that is shared by antecedents and succedents of induction axioms. This worry is for
the most part unsubstantiated, as there is a crucial difference in that the oracle is arbitrary but \emph{fixed},
whereas parameters of a scheme are universally quantified, and as such represent all numbers in the domain even in the
context of a single statement. Nevertheless, we will see that the idea that an oracle can simulate parameters works out
in certain situations, and some of our relativized separation results rely on it.

Perhaps the best way to argue that relativized separations are useful is that they show unprovability of inclusions or
reductions between rules by means of the techniques we employed elsewhere in this
paper, as all positive results we proved earlier \emph{do} relativize. This is easy to observe\footnote{The one possible
exception is that we used a couple of times the fact that every bounded sentence is provable or refutable in the base
theory. This is not literally true in the relativized setting, but it may be replaced by the weaker property that
every bounded sentence is equivalent to a Boolean combination of sentences of the form $\alpha(k)$ for standard
constants~$k$.} for the results in Sections \ref{sec:main-fragments}--\ref{sec:conserv}. For Section~\ref{sec:pps}, we
may relativize the proof systems by expanding the propositional language with a new unbounded fan-in connective
representing~$\alpha$, and then everything works out.

\begin{Thm}\label{prop:rel-pi+1indf-si}
$\pih1(\alpha)\text-\indf\nvdash T^0_2(\alpha)$, and $\pih{i+1}(\alpha)\text-\indf\nvdash S^i_2(\alpha)$ for $i>0$.
\end{Thm}
\begin{Pf}
If we fix an oracle $A\sset\omega$, then $\pih{i+1}(\alpha)\text-\indf$ follows from the set of all $\forall\sih
i(\alpha)$ sentences true in $\p{\stm,A}$. The same argument as in the proof of Theorem~\ref{prop:pi+1indf-si} then shows
that if $\pih{i+1}(\alpha)\text-\indf\vdash S^i_2(\alpha)$, then the relativized polynomial hierarchy $\ph^A$
collapses. However, it is well known that we can find~$A$ such that this does not happen \cite{yao:ph,has:olb}.

Similarly, $\pih1(\alpha)\text-\indf\vdash T^0_2(\alpha)$ implies $\ptime^A=(\tc)^A$ for every $A\sset\omega$. The proper
notion of relativized~$\tc$ corresponding to $\forall\sih1(\alpha)$-witnessing of universal extensions of~$\btc$ is
explained in Aehlig, Cook, and Nguyen~\cite{acn}, where they also exhibit an oracle separating $\cxt{NL}^A$ (hence $(\tc)^A$) from~$\ptime^A$.
\end{Pf}

\begin{Thm}\label{prop:pi+1indr-pipindf}
$\btc(\alpha)+\pih{i+1}(\alpha)\text-\indr\nvdash\pih i(\alpha)\text-\pindf$ for $i>0$, and
$\btc(\alpha)+\pih1(\alpha)\text-\indr\nvdash\sih0(\alpha)\text-\indf$.
\end{Thm}
\begin{Pf}
Assume for contradiction that $\btc(\alpha)+\pih{i+1}(\alpha)\text-\indr\vdash\pih i(\alpha)\text-\pindf$, where $i>0$.
We will argue that parameters of the $\pind$ scheme can be encoded into the oracle.

Given a term $t(x)$, let us fix a proof~$\pi$ of \pind\ for the parameter-free $\pih i(\alpha)$ formula 
\begin{equation}\label{eq:3}
\forall x_1\le t(x)\,\exists x_1\le t(x)\,\cdots Qx_i\le t(x)\,\alpha(\p{x,x_1,\dots,x_i})
\end{equation}
in $\btc(\alpha)+\pih{i+1}(\alpha)\text-\indr$, and let $\fii(x,y)$ be a $\pih
i(\alpha)$ formula of the form
\begin{equation}\label{eq:4}
\forall x_1\le t(x)\,\exists x_1\le t(x)\,\cdots Qx_i\le t(x)\,\tet(x,y,x_1,\dots,x_i),
\end{equation}
where $\tet\in\sih0(\alpha)$. We may assume without loss of generality that $y$ does not occur in~$\pi$.
If we substitute $\tet\bigl((z)_0,y,(z)_1,\dots,(z)_i\bigr)$ for $\alpha(z)$ everywhere in the proof, the result is
still a valid $\btc(\alpha)+\pih{i+1}(\alpha)\text-\indr$ proof as \indr\ allows parameters, hence the theory proves
\pind\ for $\fii(x,y)$.

This is not yet a general instance of $\pih i(\alpha)$-\pind, as all quantifiers in~$\fii$ have to be bounded by a
term in the induction variable. However, this restriction is immaterial: if $\fii(x,y)\in\pih i(\alpha)$ is arbitrary,
\pind\ for~$\fii$ follows from \pind\ for the formula $\lh x<\lh y\lor\fii(\tdive x{\lh y},y)$, which may be
equivalently rewritten so that all quantifiers are bounded by a term in $x$ alone.

Thus, $\btc(\alpha)+\pih{i+1}(\alpha)\text-\indr\vdash S^i_2(\alpha)$, but this contradicts
Theorem~\ref{prop:rel-pi+1indf-si}.

Likewise, $\btc(\alpha)+\pih1(\alpha)\text-\indr\vdash\sih0(\alpha)\text-\indf$ would imply
$\btc(\alpha)+\pih1(\alpha)\text-\indr\vdash T^0_2(\alpha)$.
\end{Pf}

We do not have an unconditional disproof of~\eqref{item:22} in its full generality, but several partial results that
come close:
\begin{Thm}\label{thm:rel-si-piindr}
Let $i\ge0$.
\begin{enumerate}
\item\label{item:42}
If $i>0$, $S^i_2(\alpha)\nvdash\btc(\alpha)+\sih i(\alpha)\text-\indr=\btc(\alpha)+\pih{i+1}(\alpha)\text-\pindr$.
\item\label{item:43}
$S^2_2(\alpha)\nvdash\btc(\alpha)+\pih2(\alpha)\text-\indr$.
\item\label{item:19}
$S^i_2(\alpha)\nvdash\pih i(\alpha)\text-\indf$.
\item\label{item:44}
$\pih i(\alpha)\text-\indr\nleq S^i_2(\alpha)$.
\end{enumerate}
\end{Thm}
\begin{Pf}
\ref{item:42}: In view of Corollary~\ref{cor:consindf}, the claim is equivalent to the fact that $T^i_2(\alpha)$ is not
$\forall\sih i(\alpha)$-conservative over~$S^i_2(\alpha)$ due to Buss and Kraj\'\i\v cek~\cite{buss-kra:sep}.

\ref{item:43}: This amounts to the $\forall\sih1(\alpha)$-non-conservativity of $T^2_2(\alpha)$ over~$S^2_2(\alpha)$,
proved by Chiari and Kraj\'\i\v cek~\cite{chi-kra:t22} (see also~\cite{chi-kra:lifting}).

\ref{item:19}: Assume that $S^i_2(\alpha)\vdash\pih i(\alpha)\text-\indf$; we will argue that $S^i_2(\alpha)\vdash\pih
i(\alpha)\text-\ind$, contradicting $S^i_2(\alpha)\ne T^i_2(\alpha)$. As in the proof of
Theorem~\ref{prop:pi+1indr-pipindf}, if $\fii(x,y)$ is a formula of the form~\eqref{eq:4}, we construct a proof of
$\fii$-\ind\ in~$S^i_2(\alpha)$ by taking a proof (not containing~$y$) of $\ind$ for the formula~\eqref{eq:3}, and
substituting $\tet\bigl((z)_0,y,(z)_1,\dots,(z)_i\bigr)$ for $\alpha(z)$. If $\fii(x,y)$ is an arbitrary $\pih
i$~formula, then $\fii$-\ind\ (with $x$ being the induction variable, and $y$ a parameter) follows from $\ind$ for the
formula $x<y\lor\fii(x-y,y)$, which is equivalent to a formula of the form~\eqref{eq:4}.

\ref{item:44} follows from \ref{item:19} using the deduction theorem.
\end{Pf}

\begin{Rem}\label{rem:rel-si-piindr}
By inspection of critical pairs of $P_R$ and~$P_T$, the net effect of
Theorems \ref{prop:rel-pi+1indf-si}, \ref{prop:pi+1indr-pipindf}, and~\ref{thm:rel-si-piindr} is that in the relativized setting:
\begin{itemize}
\item
all valid reductions between the rules from Definition~\ref{def:main} follow from Theorem~\ref{prop:basicred};
\item
all valid inclusions between theories generated by these rules follow from
Theorem~\ref{prop:basicred} and Corollary~\ref{cor:consrules}, except possibly 
\begin{align}
\label{eq:21}\btc(\alpha)&\vdash\btc(\alpha)+\sih0(\alpha)\text-\indr,\\
\intertext{or}
\label{eq:20}T&\vdash\btc(\alpha)+\pih i(\alpha)\text-\indr
\end{align}
for some~$i\ge1$, $i\ne2$, and a theory $T$ between
$\btc(\alpha)+\sih{i-1}(\alpha)\text-\indr$ and~$S^i_2(\alpha)$ (apart from the two indicated, these are
$\sih{i-1}(\alpha)$-\indf, $T^{i-1}_2(\alpha)$, $\pih i$-\pindf, and $\sih i$-\pindf).
\end{itemize}
Note that for any given~$i$, \eqref{eq:20} holds either for all the theories~$T$, or for none of them; that is, the
following are equivalent:
\begin{enumerate}
\item $S^i_2(\alpha)\vdash\btc(\alpha)+\pih i(\alpha)\text-\indr$,
\item $\btc(\alpha)+\sih{i-1}(\alpha)\text-\indr=\btc(\alpha)+\pih i(\alpha)\text-\indr$,
\item $T^i_2(\alpha)$ is $\forall\sih{i-1}(\alpha)$-conservative over~$S^i_2(\alpha)$ (or equivalently,
over~$T^{i-1}_2(\alpha)$).
\end{enumerate}
Moreover, Chiari and Kraj\'\i\v cek~\cite{chi-kra:lifting} proved that for~$i>2$, the following is also equivalent to
the above:
\begin{enumerate}
\setcounter{enumi}3
\item $T_2(\alpha)$ is $\forall\sih{i-1}(\alpha)$-conservative over~$S^i_2(\alpha)$ (or over~$T^{i-1}_2(\alpha)$).
\end{enumerate}
Likewise, \eqref{eq:21} is equivalent to the $\forall\sih0(\alpha)$-conservativity of $T^0_2(\alpha)$ over
$\btc(\alpha)$.

Even though it is commonly believed that $T^i_2(\alpha)$ is not $\forall\sih0(\alpha)$-conservative over
$S^i_2(\alpha)$ for any~$i\ge0$, it is a major open problem to improve the above-quoted results
of~\cite{buss-kra:sep,chi-kra:t22} even just by one level, thus \eqref{eq:21} and~\eqref{eq:20} are open.

In this connection, we mention a possibly interesting consequence of Theorem \ref{thm:rel-si-piindr} \ref{item:44}:
\end{Rem}
\begin{Cor}\label{cor:non-cons}
For any $i\ge1$, there is a $\forall\sih i(\alpha)$~sentence $\fii$ such that $T^i_2(\alpha)+\fii$ is not
$\forall\sih{i-1}(\alpha)$-conservative over~$S^i_2(\alpha)+\fii$.
\noproof\end{Cor}

\section{Conclusion}\label{sec:conclusion}

We have undertaken a comprehensive investigation of parameter-free and inference-rule variants of the $\sih i$
and~$\pih i$ induction and polynomial induction axioms. We found which rules and axioms reduce to other rules, and
which do not. We have seen conservation results among the systems; in particular, each of our theories can be
characterized as the $\Gamma$-fragment of some $S^i_2$ for a suitable class of sentences~$\Gamma$. We also found
equivalent expressions for our axioms and rules in terms of reflection principles for axiomatic extensions of the
quantified propositional calculi~$\G_i$, and we proved a few other results, in particular concerning nesting depth of
rules.

In some respects, the properties of our systems resemble the situation of strong theories of arithmetic $I\Sigma_n^-$
and~$I\Pi_n^-$: the $\pih i$~schemes and rules are weaker than their $\sih i$~counterparts, there are conservation
results connecting the systems to the usual theories~$S^i_2$, the parameter-free schemes do not seem to be finitely
axiomatizable, and our systems correspond to reflection principles and rules (albeit of different nature) of similar
overall shape as for the strong systems.

On the other hand, there are also notable differences. Most importantly, the hierarchies fit together in different
ways: $I\Pi_{n+1}^-$ is equiconsistent with (and $\bool(\Sigma_{n+1})$-conservative over) $I\Sigma_n^-$ and
$I\Sigma_n$, whereas in our case, $\pih i$-\ppindf\ is $\mbool(\exists\pih i\cup\forall\pih i)$-conservative under $\sih
i$-\ppindf\ and $\sih i$-\ppind. On a related note, the systems $I\Pi_{n+1}^-$ and~$I\Sigma_n$ on the same level of the
hierarchy are incomparable, and their join $I\Pi_{n+1}^-+I\Sigma_n$ has strictly stronger consistency strength---it
proves the consistency of~$I\Sigma_n$ (cf.~\cite{bekl:parfree}); no such phenomenon is possible in our setup, as all
the systems on each level of our hierarchy are included in the largest one among them, namely~$S^i_2$.

Analogously to $I\Sigma_n^-$ and $I\Pi_n^-$, it seems likely that our theories $\sih i$-\ppindf\ and $\pih i$-\ppindf\ are
not finitely axiomatizable, but we do not have any evidence for this (Question~\ref{que:finax}).
Another problem that we left open is if $T^i_2$ is $\exists\forall\pih i$-conservative over $\pih
i$-\indf, and similarly for $S^i_2$ and $\pih i$-\pindf\ (Question~\ref{que:conservp-unbd}); it would be also desirable to
prove unrelativized separation of $\pih i$-\pindf\ from $\btc+\pih{i+1}\text-\indr$ (Question~\ref{que:unrel}) under plausible
assumptions.

We tried our best to conduct an in-depth examination of parameter-free and inference-rule versions of the \ind\ and
\pind\ schemes, that also applies, by the results of Section~\ref{sec:variants}, to their common variants like \lind\
and minimization schemes. However, we left out other schemes of interest in bounded arithmetic: in particular, the
choice (aka replacement or bounded collection) scheme $\BB$ (which was studied in~\cite{cffmlm}), and analogues of
\lind\ with induction up to bounds given by more general classes of terms (including $\M{LLIND}$, etc.). Related to
$\BB$, we might be interested in variants of \ppind\ and other schemes for the non-strict $\sig i$ and~$\pii i$ formula
classes: it is well known that with parameters, the strict and non-strict \ppind\ schemes are equivalent---both define
the familiar theories $S^i_2$ and~$T^i_2$. It is however likely that the situation will get more complicated without
parameters. We also left out various combinations of our base systems such as $S^i_2+\pih i\text-\indf+\sih
i\text-\indr$.

The reason we decided not to discuss any of these potentially interesting topics
is sheer complexity: we have 10 systems per each level of the hierarchy as is, which already leads to a complex
network of relations among them. If we added more schemes and rules to the mix, the number of combinations would
multiply, rendering the global picture unmanageable. That is to say, there are certainly many aspects of these systems
that are worth further investigation, but we deem them out of scope of this paper.

\section*{Acknowledgements}
I would like to thank Andr\'es Cord\'on-Franco and F\'elix Lara-Mart\'\i n for stimulating discussions of the topic,
and for drawing my attention to~\cite{kaye:axqc}. I am grateful to the anonymous reviewer for many helpful suggestions.

\bibliographystyle{mybib}
\bibliography{bdrules}

\ifx\url\undefined {\catcode`\/=13
  \gdef/{\string/\futurelet\nexttoken\finishslash}
  \gdef\finishslash{\ifx\nexttoken/\else\penalty\relpenalty\fi}}
  \def\url{\begingroup\catcode`\~=12 \catcode`\/=13 \finishurl}
  \def\finishurl#1{\texttt{#1}\endgroup} \fi
\providecommand{\bysame}{\leavevmode\hbox to5em{\hrulefill}\thinspace}
\providecommand\bibliographyhook{}
\begin{thebibliography}{10}
\bibliographyhook

\bibitem{adam-big:isig1-}
Zofia Adamowicz and Teresa Bigorajska, \emph{Functions provably total
  in~{$I^-\Sigma_1$}}, Fundamenta Mathematicae 132 (1989), pp.~189--194.

\bibitem{acn}
Klaus Aehlig, Stephen Cook, and Phuong Nguyen, \emph{Relativizing small
  complexity classes and their theories}, Computational Complexity 25 (2016),
  no.~1, pp.~177--215.

\bibitem{bekl:indru}
Lev~D. Beklemishev, \emph{Induction rules, reflection principles, and provably
  recursive functions}, Annals of Pure and Applied Logic 85 (1997), no.~3,
  pp.~193--242.

\bibitem{bekl:parfree}
\bysame, \emph{Parameter free induction and provably total computable
  functions}, Theoretical Computer Science 224 (1999), pp.~13--33.

\bibitem{bigo:ipi1-}
Teresa Bigorajska, \emph{On {$\Sigma_1$}-definable functions provably total in
  {$I\Pi_1^-$}}, Mathematical Logic Quarterly 41 (1995), pp.~135--137.

\bibitem{bloch}
Stephen~Austin Bloch, \emph{Divide and conquer in parallel complexity and proof
  theory}, Ph.D. thesis, University of California, San Diego, 1992.

\bibitem{buss}
Samuel~R. Buss, \emph{Bounded arithmetic}, Bibliopolis, Naples, 1986, revision
  of 1985 Princeton University Ph.D. thesis.

\bibitem{buss:coll}
\bysame, \emph{Relating the bounded arithmetic and polynomial time
  hierarchies}, Annals of Pure and Applied Logic 75 (1995), no.~1--2,
  pp.~67--77.

\bibitem{bkz}
Samuel~R. Buss, Leszek~A. Ko{\l}odziejczyk, and Konrad Zdanowski,
  \emph{Collapsing modular counting in bounded arithmetic and constant depth
  propositional proofs}, Transactions of the American Mathematical Society 367
  (2015), no.~11, pp.~7517--7563.

\bibitem{buss-kra:sep}
Samuel~R. Buss and Jan Kraj{\'\i}{\v c}ek, \emph{An application of boolean
  complexity to separation problems in bounded arithmetic}, Proceedings of the
  London Mathematical Society 69 (1994), no.~3, pp.~1--21.

\bibitem{cha-kad}
Richard Chang and Jim Kadin, \emph{The {B}oolean hierarchy and the polynomial
  hierarchy: a closer connection}, SIAM Journal on Computing 25 (1996), no.~2,
  pp.~340--354.

\bibitem{chi-kra:t22}
Mario Chiari and Jan Kraj{\'\i}{\v c}ek, \emph{Witnessing functions in bounded
  arithmetic and search problems}, Journal of Symbolic Logic 63 (1998), no.~3,
  pp.~1095--1115.

\bibitem{chi-kra:lifting}
\bysame, \emph{Lifting independence results in bounded arithmetic}, Archive for
  Mathematical Logic 38 (1999), no.~2, pp.~123--138.

\bibitem{cl-tak:tc0}
Peter Clote and Gaisi Takeuti, \emph{First order bounded arithmetic and small
  boolean circuit complexity classes}, in: Feasible Mathematics {II} (P.~Clote
  and J.~B. Remmel, eds.), Progress in Computer Science and Applied Logic
  vol.~13, Birkh{\"a}user, 1995, pp.~154--218.

\bibitem{cook}
Stephen~A. Cook, \emph{Feasibly constructive proofs and the propositional
  calculus}, in: Proceedings of the 7th {A}nnual {ACM} {S}ymposium on {T}heory
  of {C}omputing, 1975, pp.~83--97.

\bibitem{cookkra}
Stephen~A. Cook and Jan Kraj{\'\i}{\v c}ek, \emph{Consequences of the
  provability of\/ {$\mathbf{NP}\subseteq\mathbf P/\mathbf{poly}$}}, Journal of
  Symbolic Logic 72 (2007), no.~4, pp.~1353--1371.

\bibitem{cook-ngu}
Stephen~A. Cook and Phuong Nguyen, \emph{Logical foundations of proof
  complexity}, Perspectives in Logic, Cambridge University Press, New York,
  2010.

\bibitem{cook-ngu:err}
\bysame, \emph{Corrections for \cite{cook-ngu}},
  \url{http://www.cs.toronto.edu/~sacook/homepage/corrections.pdf}, 2013.

\bibitem{cffmlm}
Andr{\'e}s Cord{\'o}n-Franco, Alejandro Fern{\'a}ndez-Margarit, and
  Francisco~F{\'e}lix Lara-Mart{\'\i}n, \emph{Existentially closed models and
  conservation results in bounded arithmetic}, Journal of Logic and Computation
  19 (2009), no.~1, pp.~123--143.

\bibitem{cf-lm:loc-ind}
Andr{\'e}s Cord{\'o}n-Franco and Francisco~F{\'e}lix Lara-Mart{\'\i}n,
  \emph{Local induction and provably total computable functions}, Annals of
  Pure and Applied Logic 165 (2014), no.~9, pp.~1429--1444.

\bibitem{hp}
Petr H{\'a}jek and Pavel Pudl{\'a}k, \emph{Metamathematics of first-order
  arithmetic}, Perspectives in Mathematical Logic, Springer, 1993, second
  edition 1998.

\bibitem{has:olb}
Johan H{\aa}stad, \emph{Almost optimal lower bounds for small depth circuits},
  in: Randomness and Computation (S.~Micali, ed.), Advances in Computing
  Research: A Research Annual vol.~5, JAI Press, 1989, pp.~143--170.

\bibitem{ej:hash}
Emil Je{\v r}{\'a}bek, \emph{Approximate counting by hashing in bounded
  arithmetic}, Journal of Symbolic Logic 74 (2009), no.~3, pp.~829--860.

\bibitem{ej:vnc}
\bysame, \emph{On theories of bounded arithmetic for {$\mathit{NC}^1$}}, Annals
  of Pure and Applied Logic 162 (2011), no.~4, pp.~322--340.

\bibitem{ej-ngu:strict}
Emil Je{\v r}{\'a}bek and Phuong Nguyen, \emph{Simulating non-prenex cuts in
  quantified propositional calculus}, Mathematical Logic Quarterly 57 (2011),
  no.~5, pp.~524--532.

\bibitem{joh-pol:d1cr}
Jan Johannsen and Chris Pollett, \emph{On the {$\Delta^b_1$}-bit-comprehension
  rule}, in: Logic {C}olloquium '98: Proceedings of the 1998 {ASL} {E}uropean
  {S}ummer {M}eeting held in {P}rague, {C}zech {R}epublic (S.~R. Buss,
  P.~H{\'a}jek, and P.~Pudl{\'a}k, eds.), ASL, 2000, pp.~262--280.

\bibitem{kaye:parf}
Richard Kaye, \emph{Parameter free induction in arithmetic}, in: Proceedings of
  the 5th {E}aster {C}onference on {M}odel {T}heory, Sektion Mathematik der
  Humboldt-Universit{\"a}t zu Berlin, 1987, pp.~70--81, {S}eminarbericht
  {N}r.~93.

\bibitem{kaye:axqc}
\bysame, \emph{Axiomatizations and quantifier complexity}, in: Proceedings of
  the 6th {E}aster {C}onference on {M}odel {T}heory, Sektion Mathematik der
  Humboldt-Universit{\"a}t zu Berlin, 1988, pp.~65--84, {S}eminarbericht
  {N}r.~98.

\bibitem{kaye:dio}
\bysame, \emph{Diophantine induction}, Annals of Pure and Applied Logic 46
  (1990), no.~1, pp.~1--40.

\bibitem{kpd}
Richard Kaye, Jeff Paris, and Costas Dimitracopoulos, \emph{On parameter free
  induction schemas}, Journal of Symbolic Logic 53 (1988), no.~4,
  pp.~1082--1097.

\bibitem{kra:witlog}
Jan Kraj{\'\i}{\v c}ek, \emph{Fragments of bounded arithmetic and bounded query
  classes}, Transactions of the American Mathematical Society 338 (1993),
  no.~2, pp.~587--598.

\bibitem{book}
\bysame, \emph{Bounded arithmetic, propositional logic, and complexity theory},
  Encyclopedia of Mathematics and Its Applications vol.~60, Cambridge
  University Press, 1995.

\bibitem{krpu}
Jan Kraj{\'\i}{\v c}ek and Pavel Pudl{\'a}k, \emph{Quantified propositional
  calculi and fragments of bounded arithmetic}, Zeit\-schrift f{\"u}r
  ma\-the\-ma\-ti\-sche Lo\-gik und Grund\-la\-gen der Ma\-the\-ma\-tik 36
  (1990), no.~1, pp.~29--46.

\bibitem{kpt}
Jan Kraj{\'\i}{\v c}ek, Pavel Pudl{\'a}k, and Gaisi Takeuti, \emph{Bounded
  arithmetic and the polynomial hierarchy}, Annals of Pure and Applied Logic 52
  (1991), no.~1--2, pp.~143--153.

\bibitem{par-wil:counting}
Jeff~B. Paris and Alex~J. Wilkie, \emph{Counting problems in bounded
  arithmetic}, in: Methods in Mathematical Logic (C.~A. Di~Prisco, ed.),
  Lecture Notes in Mathematics vol. 1130, Springer, 1985, pp.~317--340.

\bibitem{sk-th:gi}
Alan Skelley and Neil Thapen, \emph{The provably total search problems of
  bounded arithmetic}, Proceedings of the London Mathematical Society 103
  (2011), no.~1, pp.~106--138.

\bibitem{yao:ph}
Andrew C.-C. Yao, \emph{Separating the polynomial-time hierarchy by oracles},
  in: Proceedings of the 26th {A}nnual {IEEE} {S}ymposium on {F}oundations of
  {C}omputer {S}cience (R.~E. Tarjan, ed.), 1985, pp.~1--10.

\bibitem{zamb:notes}
Domenico Zambella, \emph{Notes on polynomially bounded arithmetic}, Journal of
  Symbolic Logic 61 (1996), no.~3, pp.~942--966.

\end{thebibliography}
\end{document}